\documentclass[12pt]{amsart}
\usepackage{amscd,amssymb,latexsym}
\usepackage{tikz}
\usepackage{stmaryrd}
\usepackage{fullpage}
\usepackage{hyperref}
\usepackage{resizegather}
\usepackage{caption}
\usepackage{color}
\usepackage{color, colortbl}
\usepackage{array}
\def\mystrut{\rule[0ex]{0ex}{2.15ex}}

\usetikzlibrary{shapes.geometric}
\input{xypic}
\newcommand{\Mdef}[2]{\newcommand{#1}{\relax \ifmmode #2 \else $#2$\fi}}

\newcommand{\sm}{\wedge}

\newcommand{\rost}{\bigstar}

\newcommand{\Hom}{\mathrm{Hom}}

\newcommand{\Ext}{\mathrm{Ext}}
\Mdef{\bhom}{\mathbf{\hat{H}om}}

\Mdef{\Mod}{\mathrm{mod}}

\newcommand{\st}{\; | \;}

\numberwithin{equation}{section}

\theoremstyle{definition}

\newcommand{\qqed}{\qed \\[1ex]}

\Mdef{\PH} {\Phi^H}
\Mdef{\PK} {\Phi^K}
\Mdef{\PL} {\Phi^L}
\Mdef{\PT} {\Phi^{\T}}

\Mdef{\ef}{E{\cF}_+}
\Mdef{\etf}{\widetilde{E}{\cF}}
\Mdef{\eg}{E{G}_+}
\Mdef{\etg}{\tilde{E}{G}}
\newcommand{\Mpi}{\underline{\pi}}

\Mdef{\infl}{\mathrm{inf}}
\Mdef{\defl}{\mathrm{def}}
\Mdef{\res}{\mathrm{res}}
\Mdef{\ind}{\mathrm{ind}}
\Mdef{\coind}{\mathrm{coind}}

\Mdef{\univ}{\mathcal{U}}

\Mdef{\Fp}{\mathbb{F}_p}
\Mdef{\Zpinfty}{\Z /p^{\infty}}
\Mdef{\Zpadic}{\Z_p^{\wedge}}

\newcommand{\bi}{\begin{itemize}}
\newcommand{\be}{\begin{enumerate}}
\newcommand{\bc}{\begin{center}}
\newcommand{\bd}{\begin{description}}
\newcommand{\ei}{\end{itemize}}
\newcommand{\ee}{\end{enumerate}}
\newcommand{\ec}{\end{center}}
\newcommand{\ed}{\end{description}}

\newcommand{\lra}{\longrightarrow}
\newcommand{\lla}{\longleftarrow}

\Mdef{\we}{\mathbf{we}}
\Mdef{\fib}{\mathbf{fib}}
\Mdef{\cof}{\mathbf{cof}}
\Mdef{\BI}{\mathcal{BI}}

\Mdef{\A}{\mathbb{A}}
\Mdef{\B}{\mathbb{B}}
\Mdef{\C}{\mathbb{C}}
\Mdef{\D}{\mathbb{D}}
\Mdef{\E}{\mathbb{E}}
\Mdef{\T}{\mathbb{T}}
\Mdef{\F}{\mathbb{F}}
\Mdef{\G}{\mathbb{G}}
\Mdef{\I}{\mathbb{I}}
\Mdef{\N}{\mathbb{N}}
\Mdef{\Q}{\mathbb{Q}}
\Mdef{\R}{\mathbb{R}}
\Mdef{\bbS}{\mathbb{S}}
\Mdef{\Z}{\mathbb{Z}}

\Mdef{\bA}{\mathbb{A}}
\Mdef{\bB}{\mathbb{B}}
\Mdef{\bC}{\mathbb{C}}
\Mdef{\bD}{\mathbb{D}}
\Mdef{\bE}{\mathbb{E}}
\Mdef{\bF}{\mathbb{F}}
\Mdef{\bG}{\mathbb{G}}
\Mdef{\bH}{\mathbb{H}}
\Mdef{\bI}{\mathbb{I}}
\Mdef{\bJ}{\mathbb{J}}
\Mdef{\bK}{\mathbb{K}}
\Mdef{\bL}{\mathbb{L}}
\Mdef{\bM}{\mathbb{M}}
\Mdef{\bN}{\mathbb{N}}
\Mdef{\bO}{\mathbb{O}}
\Mdef{\bP}{\mathbb{P}}
\Mdef{\bQ}{\mathbb{Q}}
\Mdef{\bR}{\mathbb{R}}
\Mdef{\bS}{\mathbb{S}}
\Mdef{\bT}{\mathbb{T}}
\Mdef{\bU}{\mathbb{U}}
\Mdef{\bV}{\mathbb{V}}
\Mdef{\bW}{\mathbb{W}}
\Mdef{\bX}{\mathbb{X}}
\Mdef{\bY}{\mathbb{Y}}
\Mdef{\bZ}{\mathbb{Z}}

\Mdef{\cA}{\mathcal{A}}
\Mdef{\cB}{\mathcal{B}}
\Mdef{\cC}{\mathcal{C}}
\Mdef{\mcD}{\mathcal{D}} % Something funny about \cD.
\Mdef{\cE}{\mathcal{E}}
\Mdef{\cF}{\mathcal{F}}
\Mdef{\cG}{\mathcal{G}}
\Mdef{\mcH}{\mathcal{H}} % There's something funny about \cH: it
\Mdef{\cI}{\mathcal{I}}
\Mdef{\cJ}{\mathcal{J}}
\Mdef{\cK}{\mathcal{K}}
\Mdef{\mcL}{\mathcal{L}}% There's something funny about \cL: it
\Mdef{\cM}{\mathcal{M}}
\Mdef{\cN}{\mathcal{N}}
\Mdef{\cO}{\mathcal{O}}
\Mdef{\cP}{\mathcal{P}}
\Mdef{\cQ}{\mathcal{Q}}
\Mdef{\mcR}{\mathcal{R}}% There's something funny about \cR: it
\Mdef{\cS}{\mathcal{S}}
\Mdef{\cT}{\mathcal{T}}
\Mdef{\cU}{\mathcal{U}}
\Mdef{\cV}{\mathcal{V}}
\Mdef{\cW}{\mathcal{W}}
\Mdef{\cX}{\mathcal{X}}
\Mdef{\cY}{\mathcal{Y}}
\Mdef{\cZ}{\mathcal{Z}}

\Mdef{\ca}{\mathcal{a}}
\Mdef{\ct}{\mathcal{t}}

\Mdef{\At}{\tilde{A}}
\Mdef{\Bt}{\tilde{B}}
\Mdef{\Ct}{\tilde{C}}
\Mdef{\Et}{\tilde{E}}
\Mdef{\Ht}{\tilde{H}}
\Mdef{\Kt}{\tilde{K}}
\Mdef{\Lt}{\tilde{L}}
\Mdef{\Mt}{\tilde{M}}
\Mdef{\Nt}{\tilde{N}}
\Mdef{\Pt}{\tilde{P}}

\Mdef{\tA}{\tilde{A}}
\Mdef{\tB}{\tilde{B}}
\Mdef{\tC}{\tilde{C}}
\Mdef{\tE}{\tilde{E}}
\Mdef{\tH}{\tilde{H}}
\Mdef{\tK}{\tilde{K}}
\Mdef{\tL}{\tilde{L}}
\Mdef{\tM}{\tilde{M}}
\Mdef{\tN}{\tilde{N}}
\Mdef{\tP}{\tilde{P}}

\Mdef{\ft}{\tilde{f}}
\Mdef{\xt}{\tilde{x}}
\Mdef{\yt}{\tilde{y}}

\Mdef{\Ab}{\overline{A}}
\Mdef{\Bb}{\overline{B}}
\Mdef{\Cb}{\overline{C}}
\Mdef{\Db}{\overline{D}}
\Mdef{\Eb}{\overline{E}}
\Mdef{\Fb}{\overline{F}}
\Mdef{\Gb}{\overline{G}}
\Mdef{\Hb}{\overline{H}}
\Mdef{\Ib}{\overline{I}}
\Mdef{\Jb}{\overline{J}}
\Mdef{\Kb}{\overline{K}}
\Mdef{\Lb}{\overline{L}}
\Mdef{\Mb}{\overline{M}}
\Mdef{\Nb}{\overline{N}}
\Mdef{\Ob}{\overline{O}}
\Mdef{\Pb}{\overline{P}}
\Mdef{\Qb}{\overline{Q}}
\Mdef{\Rb}{\overline{R}}
\Mdef{\Sb}{\overline{S}}
\Mdef{\Tb}{\overline{T}}
\Mdef{\Ub}{\overline{U}}
\Mdef{\Vb}{\overline{V}}
\Mdef{\Wb}{\overline{W}}
\Mdef{\Xb}{\overline{X}}
\Mdef{\Yb}{\overline{Y}}
\Mdef{\Zb}{\overline{Z}}

\Mdef{\db}{\overline{d}}
\Mdef{\hb}{\overline{h}}
\Mdef{\qb}{\overline{q}}
\Mdef{\rb}{\overline{r}}
\Mdef{\tb}{\overline{t}}
\Mdef{\ub}{\overline{u}}
\Mdef{\vb}{\overline{v}}

\Mdef{\hc}{\hat{c}}
\Mdef{\he}{\hat{e}}
\Mdef{\hf}{\hat{f}}
\Mdef{\hA}{\hat{A}}
\Mdef{\hH}{\hat{H}}
\Mdef{\hJ}{\hat{J}}
\Mdef{\hM}{\hat{M}}
\Mdef{\hP}{\hat{P}}
\Mdef{\hQ}{\hat{Q}}

\Mdef{\thetab}{\overline{\theta}}
\Mdef{\phib}{\overline{\phi}}

\Mdef{\uA}{\underline{A}}
\Mdef{\uB}{\underline{B}}
\Mdef{\uC}{\underline{C}}
\Mdef{\uD}{\underline{D}}
\Mdef{\uvb}{\underline{\vb}}
\Mdef{\ul}{\underline{l}}

\Mdef{\bolda}{\mathbf{a}}
\Mdef{\boldb}{\mathbf{b}}
\Mdef{\bfD}{\mathbf{D}}

\Mdef{\fm}{\frak{m}}
\Mdef{\fp}{\frak{p}}

\Mdef{\eps}{\epsilon}

\newcommand{\Ftwo}{\mathbb{F}_2}

\newcommand{\BPRn}{BP\R \langle n \rangle}

\newcommand{\cell}{\mathrm{Cell}}
\usepackage{color}
\definecolor{todo}{rgb}{1,0,0}

\definecolor{Gray}{gray}{0.85}
\definecolor{LightCyan}{rgb}{0.88,1,1}

\begin{document}

\baselineskip 15.00pt

\title{The representation-ring-graded local cohomology spectral sequence for $BP\R\langle 3 \rangle$}

\author{J.P.C. Greenlees}
\address{School of Mathematics and Statistics, Hicks Building, Sheffield S3 7RH. UK.}
\email{j.greenlees@sheffield.ac.uk}

\author{Dae-Woong Lee}
\address{Department of Mathematics, and Institute of Pure and Applied Mathematics, Chonbuk National University, 567 Baekje-daero, Deokjin-gu, Jeonju-si, Jeollabuk-do 54896, Republic of Korea}
\email{dwlee@jbnu.ac.kr}
\date{15 August 2017}


\begin{abstract}
The purpose of this paper is to examine the calculational consequences
of the duality proved by the first author and Meier in \cite{BPRn},
making them explicit in one more example.
We give an explicit description of the behaviour of the spectral
sequence in the title, with pictures. 
\end{abstract}

\thanks{
The authors are grateful to L.Meier for comments on a draft of this paper. 
The second author was supported by Basic Science Research Program through the
National Research Foundation of Korea (NRF) funded by the Ministry of
Education (NRF-2015R1D1A1A09057449).
}
\maketitle

%\tableofcontents

\newcommand{\uk}{\underline{k}}
\section{Preamble}
\subsection{Summary}
It was shown by the first author and Meier \cite{BPRn} that Gorenstein
duality for $\BPRn$ gives rise to a certain 2-local  local cohomology
spectral sequence graded over the real representation ring. 
%With several explanations to follow, 
In this paper we will describe
the behaviour of the spectral sequence
$$H^*_J(BP\R\langle 3 \rangle^Q_{\rost})\Rightarrow
\Sigma^{-(16+9\sigma)} \pi^Q_{\rost}(\Z_{(2)}^{BP\R\langle 3 \rangle}), $$
where the homotopy of the Anderson dual $\Z_{(2)}^{BP\R\langle 3 \rangle}$ lives in a short exact
sequence
$$0\lra
\Ext_{\Z}(\pi^Q_{\rost}(\Sigma {BP\R\langle 3 \rangle}), \Z_{(2)}) \lra
\pi^Q_{\rost}(\Z_{(2)}^{BP\R\langle 3 \rangle})\lra \Hom_{\Z}(
\pi^Q_{\rost}({BP\R\langle 3 \rangle}), \Z_{(2)})\lra 0. $$

There are many things in this statement that need a proper
introduction later, but here are the absolute essentials. 
There is a well known complex orientable spectrum $BP\langle 3\rangle$ at the prime 2 (with coefficient ring
$\Z_{(2)}[v_1, v_2]$), and the spectrum 
 $BP\R\langle 3\rangle$ is a version that takes into account the action of the group
$Q$ of order $2$ by complex conjugation.  The subscript $\rost$ refers to  
grading over the real representation ring $RO(Q)=\{ x+y\sigma\st x,
y\in \Z\}$,  where $\sigma$ is the sign representation on $\R$. The
ideal $J$ is generated by elements $\vb_1, \vb_2$, and $H^*_J$ denotes
local cohomology in the sense of Grothendieck. 

We will introduce the statement and indicate its interest
by giving a quasi-historical account (see \cite{hica} for more
context).  We will make the assertions more
precise in the process.

\subsection{Evolution}
The classical local cohomology theorem \cite{KEG, groupca} with shift $a$ for a cohomology theory
$R^*(\cdot )$  gives an approach to $R_*(BG)$ through the local
cohomology of $R^*(BG)$ for a compact Lie group $G$: there is a spectral
sequence
$$H^*_J(R^*(BG))\Rightarrow \Sigma^a R_*(BG),$$
where $J=\ker (R^*(BG)\lra R^*)$.
When $R^*(\cdot)$ is ordinary cohomology with coefficients in a field
$k$, the reverse universal coefficient theorem takes a very
simple form: $$R_*(BG)=\Hom_k(R^*(BG),  k)$$ (amounting to the statement
that $R\sm BG_+$ is Brown-Comenetz dual to $F(BG_+, R)$). The Local
Cohomology Theorem (LCT for short) \cite{groupca} then shows  that
the ring $R^*(BG)$ is very special \cite{ringlct} (for example it is Gorenstein in
codimension 0).

More generally $R^*(BG)$ (which is the coefficient ring of
Borel-equivariant $R$ cohomology),  may be replaced by the equivariant
coefficient ring $R^*_G$. The natural replacement of $R_*(BG)$ is then
$R^G_*(EG)$, and there are interesting instances of the LCT of that
type. The new complication here is that when $G$ is not finite, there
may be a twist arising from the isomorphism
$$R^G_*(S^{-ad(G)}\sm EG_+)\cong R_*(BG_+),  $$
where $ad (G)$ denotes the adjoint representation. 
In particular, it becomes natural to consider not just the
$\Z$-graded coefficient ring $R_G^*$ but the $RO(G)$-graded version.

The lesser complication that we need here is that the ring $k$ we work
over will not be a field but rather a localization of the
integers. This means that $R\sm BG_+$ (as {\em Brown-Comenetz} dual \cite{BC76}) should
be replaced by the equivariant {\em Anderson} dual $k^R$. We will explain
this in more detail in Subsection \ref{subsec:GorD}.

\subsection{Equivariant cohomology theories}
We will consider  equivariant cohomology theories $R_G^*(X)$
represented by a (genuine) $G$-spectrum $R$:
$$R_G^*(X)=[X, R]^G_*. $$
We will often consider orthogonal representations $V$ of $G$, and
their one point compactifications $S^V$. The sphere $S^V$ is
invertible as a $G$-spectrum, and we will consider $RO(G)$-graded
rings $R^G_{\rost}$ (with $R^G_V=[S^V, R]^G$ in degree $V$).

\subsection{Mackey functors}
The equivariant homotopy groups $\pi^H_n(X)$  for
various subgroups $H$ are related by induction and restriction maps,
in the sense that
$$\Mpi_n^G(X)=\{ \pi^K_n(X)\}_K. $$
is a Mackey functor.

% We will use Dress’s formulation of Mackey functors for a finite group
%$ G$  \cite{Dress} as given by a covariant and a contravariant functor
%on finite $G$-sets subject to the Mackey condition.
%The category of $G$-spectra includes the orbits $G/H_+$, and we write
%$$\SO_G$ for the stable orbit category they generate.
%A Mackey functor
%$M$ for a group
%$G$ is therefore determined by its values on the transitive $G$-sets
%$G/H$ together with certain structure maps.
 %If $K \subseteq  H  $ the projection $\pi^H_K: G/K \lra G/H$ induces
 %a restriction map
%$(\pi^H_K)^* = \res^H_K : M(G/H)\lra  M(G/K)$ and an induction map
%$(\pi^H_K)_* = \ind^H_K : M (G/K )\lra  M (G/H )$,  and right%
%multiplication  $R_g : G/H \lra G/H^g$ induces an action of $W_G(H) =
%N_G(H)/H$  on $M(G/H)$.

We will need to refer to just two very simple Mackey functors.
For an abelian group $k$, the {\em constant Mackey functor} $\uk$ has the value $k$ at all subgroups, the
restriction maps are the identity and the transfer maps are
multiplication by the index. The dual Mackey functor $\uk^*$ has all
the transfer maps the identity and the restriction maps are
multiplication by the index.

\subsection{Gorenstein duality}
\label{subsec:GorD}
Suppose  $R$ is a connective $G$-equivariant ring spectrum with  0th
homotopy Mackey functor constant at the abelian group $k$:
$\Mpi_0^G(R)=\uk$. Killing homotopy groups in $R$-modules gives  a map
$R\lra H\uk$.

 We say that $R$ is {\em Gorenstein} of shift $V$ if there is an equivalence
$$\Hom_R(H\uk, R)\simeq \Sigma^V H\uk $$
of $G$-equivariant $R$-modules. If $\cell_{H\uk}R\lra R$ is the 
$G-H\uk$-cellularization of $R$ (i.e., the best approximation to $R$
built with objects $G/H_+\sm H\uk$),  this means
$$\Hom_R(H\uk, \cell_{H\uk} R)\simeq \Sigma^V H\uk .$$

 In practice we are interested in a closely
related duality statement.  This involves considering the Anderson dual
$G$-spectrum $k^R$;  
%The Anderson dual  can be defined by choosing an injective resolution of the
%abelian group $k$, but
for our purpose the essential fact about $k^R$ is that its
homotopy groups are dual to those of $R$ in the sense that, for any
$G$-spectrum $T$,  there is a short exact sequence
% \lra I\lra J\lra 0$ of abelian groups.
$$ 0\lra \Ext_{\Z} (\pi^G_*(\Sigma T\sm R), k)\lra     [T, k^R]^G_*\lra  \Hom_{\Z}
(\pi^G_*(T\sm R), k)\lra   0. $$
Note in particular the special cases where $T=G/H_+$, which means  that
$$\Hom_R(H\uk , k^R)\simeq H\uk^*, $$
where $\uk^*$ is the dual of the constant Mackey functor (the one in
which the induction maps are all the identity).

Note also the special case with $T=S^V$, showing that when the abelian
groups are projective, $\pi_V^G(k^R)=\Hom_{\Z}(\pi_{-V}^G(R), \Z)$.

Now  for many small groups $G$ the Anderson duality is a twist for
integral homology in the sense that there
is a virtual representation $\gamma$ so that
$$H\uk \simeq \Sigma^{\gamma}H\uk^*$$
(For example if $G$ is cyclic we may take
$\gamma =\epsilon-\alpha$, where $\epsilon$ is the trivial
representation on $\C$ and $\alpha $ is a faithful representation of $G$
on $\C$). We then have
$$\Hom_R(H\uk , \cell_{H\uk} R)\simeq \Sigma^V H\uk \simeq
\Sigma^{V+\gamma}H\uk^*
\simeq \Hom_R(H\uk, \Sigma^{V+\gamma} k^R). $$
Now $\Sigma^{V+\gamma}k^R$ is  
$G-H\uk$-cellular $G$-spectra (this is not hard to see using the fact
that Anderson duality is a twist for integral homology), so that under orientability hypotheses, one may hope that
$\Hom_R(H\uk, \cdot )$ can be removed from the equivalence.

We say that $R$ has {\em Gorenstein duality} of shift $W$ \cite{DGI06} if
$$\cell_{H\uk} R \simeq \Sigma^W k^R,   $$
so that given orientability and the hypothesis on Eilenberg-MacLane
spectra if $R$ is Gorenstein of shift $V$ then it also has Gorenstein
duality  of shift $W=V+\gamma$ \cite{BPRn}. The purpose of the present
paper is to study this duality in a case established in 
\cite{BPRn}.

\subsection{The local cohomology spectral sequence}
Now in favourable circumstances, one may construct the  
$G-H\uk$-cellularization by algebraic means, using the ideal
$$\hat{J}=\ker (R^G_{\rost}\lra H\uk^G_{\rost})$$
or a smaller ideal $J$ with the same radical:
$$\cell_{H\uk} M\simeq \Gamma_JM.$$
In this case, Gorenstein duality with shift $W$ gives rise to a local
cohomology spectral sequence
$$H^*_J(R^G_{\rost})\Rightarrow \Sigma^{W}\pi^G_{\rost}(k^R). $$
We are going to study this spectral sequence in a case established in 
\cite{BPRn}.

\subsection{The spectrum $\BPRn$}
We now specialize the group $G$ to be the group $Q=Gal (\C | \R)$\footnote{The use of the letter $Q$ has unaccountably
  attracted considerable attention!  It reflects the fact that if we start with a group
$\Gamma$ of equivariance,  adding in a Galois group $Q$ (as in \cite{AtiyahSegal, Karoubi}) we will 
need to consider $\hat{\Gamma}$-equivariant spectra for an extension 
$1\lra \Gamma \lra \hat{\Gamma}\lra Q\lra 1$, so that $Q$ is a
quotient. In our case $\Gamma$ is trivial, obscuring the motivation.}
of order 2. 
 We write $\sigma$ for the sign representation of $ Q$ on
$\R$, $1$ for the  trivial representation, and $\rho=1+\sigma$  for the real regular  representation. We will be grading our groups over the real
representation ring
$$RO(Q) = \{ x + y\sigma  \st  x, y \in \Z\}.$$
When we   draw pictures, they will be displayed in the plane with the $\Z$ axis  horizontal and the $\sigma$  axis vertical, so that $x+y\sigma$ gets displayed at the point with cartesian coordinates $(x,y)$.

The spectrum $\BPRn$ is a $Q$-spectrum formed from the  Real bordism
spectrum $M\R$ considered by Araki, Landweber and Hu-Kriz \cite{HK01}.
Localizing at $2$, a $Q$-equivariant refinement of the Quillen
idempotent defines the $Q$-spectrum $BP\R$. The Real spectra we are
interested in are quotients of
$BP\R$ by homotopy elements $\vb_i$ (of degree $(2^i-1)\rho$) which are
equivariant refinements of the familiar elements $v_i$ (of degree
$(2^i-1)2$). There are various possible choices of such elements, but
these  will not affect what we say.  For example, we can follow \cite{HK01} and \cite{Hu02} and define
$$\BPRn = BP\R/(\vb_{n+1},\vb_{n+2}, \ldots ).$$
For the present, the most important thing is that in degrees which are
multiples of $\rho$ the coefficients are rather familiar
$$\BPRn^Q_{*\rho}=\Z_{(2)} [\vb_1, \ldots , \vb_n];$$
furthermore the augmentation ideal is the radical of $J=(\vb_1, \ldots , \vb_n)$.   A little more history is described in \cite{BPRn}, but see
 \cite{Ati66}, \cite{L68}, 
 \cite{Dug05}, \cite{BG10}, 
% \cite{LO16},  \cite{Lor16} 
and \cite{HM}.

\subsection{The local cohomology spectral sequence for $\BPRn$}
In \cite[Theorem 1.1]{BPRn}, the first author and L. Meier showed when $G=Q$ is
the Galois group $Gal(\C|\R)$ of order 2, the $Q$-spectrum  $\BPRn$
has Gorenstein duality of shift $-[D_n \rho +n+2(1-\sigma)]$, where
$$D_n\rho=|\vb_1|+\cdots +|\vb_n|=(2^{n+1}-n-2)\rho,$$ and that
consequently there is a Local Cohomology spectral sequence
$$E^2_{*,\rost}=H^*_J(\BPRn^Q_{\rost})\Rightarrow \Sigma^{-[D_n \rho
  +n+2(1-\sigma)]}\pi^Q_{\rost}(\Z_{(2)}^{\BPRn}),  $$
where the differentials are of the standard homological form $d_i:
E^i_{n,V}\lra E^i_{n-i, V+i-1}$ changing the total $RO(Q)$-degree by
adding the integer $-1$.  The result was illustrated in  \cite{BPRn}  by working out the exact behaviour of the
LCT for $n=1$ and $n=2$.

The local cohomology spectral sequence display in standard homological
form has the $i$th local cohomology $H^i_{J}$ displayed on the
$-i$th column. If $J$ has $n$ generators, this means it only has
non-zero entries on the columns $-n, -n+1, \ldots , -1, 0$, so the
last possible differential is $d_n$.

For $\BPRn$, the relevant ideal is  $J=(\vb_1, \ldots ,
\vb_n)$. Accordingly, if $n=1$ there can be no differentials, but
there are some additive and multiplicative extensions. If $n=2$ it was
shown in \cite{BPRn} that  there  are precisely three nonzero $d_2$
differentials up to periodicity.

The present paper describes the behaviour of the spectral sequence for
$n=3$. This time there are both nontrivial $d_2$ and nontrivial $d_3$
differentials. This makes clear the
expectation that for all $n$ the differentials $d_2, d_3, \ldots ,
d_n$ will all be non-zero, and we prove $d_n\neq 0$ for all $n \geq 2$ in Section \ref{sec:diff}.  We also describe the additive and multiplicative extensions that
occur.

The other fundamental feature is that since $J$ is generated by
elements whose degree is a multiple of $\rho$, the local cohomology
can be calculated taking one diagonal at a time (i.e., fixing $\delta$
and focusing on gradings of the form $\delta +*\rho$). Since
Anderson duality also preserves diagonals, this means that the
Gorenstein duality breaks up into dualities by diagonal.
Poincar\'e duality for an $n$-manifold pairs $i$th homology and
$(n-i)$th cohomology. Similarly, Gorenstein duality behaves as if we
have a manifold of dimension $2^{n+1}-4$: it pairs the local cohomology
of the  $\delta$th diagonal of  $\BPRn^Q_{\rost}$ (for $0\leq \delta \leq 2^{n+1}$)
with the $(2^{n+1}-4-\delta -\epsilon)$th diagonal
$\BPRn^Q_{\rost}$ (where $\epsilon$ is a small homological correction, just as
 happens with torsion in integral Poincar\'e duality).
This was described for  $n=1$ (like a $4$-manifold) and $n= 2$ (like a
$12$-manifold)  in \cite{BPRn}, and we describe it
here for $n= 3$ (like a $28$-manifold).

\bigskip

\section{The ring depicted}
The $RO(Q)$-graded ring $\BPRn^Q_{\rost}$ is displayed in \cite{BPRn}
for $n=1,2$ and for $n=0,1$ in many places (but with the present
notation in  \cite{zkr}). In this section we display it  for $n=3$.

The general picture is that the homotopy fixed point groups $\BPRn^{hQ}_{\rost}$ have a periodicity
element $U$ of degree $2^{n+1}(1-\sigma)$, so that it is of the form
$$\BPRn^{hQ}_{\rost}=BB[U,U^{-1}]$$
for a certain `basic block' $BB$ (depending on $n$), which is a module over
$\Z_{(2)} [a] [\vb_1, \ldots, \vb_n]$, where $a$ is of degree $-\sigma$; see
\cite[Proposition 4.1]{BPRn} for details. Most
elements of $BB$ are $a$-power torsion (in fact the only exception is
a single copy of $\Z [a]/(2a)$). The process of recovering the actual
fixed points from the homotopy fixed points corresponds to taking
$(a)$-local cohomology (see \cite[Subsection 11.B]{BPRn}), so 
we obtain a negative block $NB$ which is the same as $BB$ except in
degrees $x+y\sigma$ with $x=0,-1$. We then find
$$\BPRn^{Q}_{\rost} =BB[U]\oplus U^{-1}\cdot NB[U^{-1}]. $$
As mentioned, $NB$ is very similar to $BB$, so $\BPRn^Q_{\rost}$ is nearly
$U$-periodic.

In  the specific case $n=3$, we find $U$ is of degree
$16(1-\sigma)$ and we may display
$$
BP\mathbb R\langle 3 \rangle^{Q}_{\rost}=BB[U] \oplus U^{-1} \cdot NB[U^{-1}].
$$
We show this in the $x+y\sigma$ plane. The details will be
described in due course, but for now note that small red dots are copies
of $\Ftwo$ and larger circles or squares are copies of $\Z$. It is
apparent that almost all entries are in the region $x+y\geq 0$. In the
left half-plane ($x<0$) this is precisely true, and in the right
half-plane ($x\geq 0$) it is true except for vertical strings hanging
down. The other obvious feature is the $U$ almost-periodicity referred
to above.

%\newpage

\begin{figure}
\centering
$$\begin{tikzpicture}[scale =0.425]
%\clip (-9, -8.6) rectangle (5, 8.5);
\draw[step=0.5, gray, very thin] (-20,-20) grid (20, 20);
\draw (-8.0,-3.0 ) node[anchor=east, draw=orange]{ \Large{ $BP\mathbb{R} \langle 3 \rangle^{Q}_{\rost}$  }  };

\foreach \y in {1,2,3,4,5,6,7,8,9,10,11,12,13,14,15,16,17,18,19,20,21,22,23,24,25,26,27,28,29,30,31,32,33,34,35,36,37,38,39,40}
\node at (0,0-\y/2) [fill=red, inner sep=1pt, shape=circle, draw] {};

\draw [->] (0,0)-- (20,20);
\node at (0,0)  [shape = rectangle, draw]{};
\draw (0,0) node[anchor=east]{1};

\draw[->][red] (0,-0.5)--(20,19.5);
\draw[->][red] (0,-1.0)--(20,19.0);

\draw[->][blue] (0,0.5-2)--(20,18.5);
\draw[->][blue] (0,-0.10-2)--(20,17.90);
\draw[->][blue] (0,-0.5-2)--(20,17.50);
\draw[->][blue] (0,-1.0-2)--(20,17.00);

\draw[->][green] (0,-1.5-2)--(20,16.50);
\draw[->][green] (0,-2.10-2)--(20,15.90);
\draw[->][green] (0,-2.60-2)--(20,15.40);
\draw[->][green] (0,-3.10-2)--(20,14.90);
\draw[->][green] (0,-3.5-2)--(20,14.50);
\draw[->][green] (0,-4.10-2)--(20,13.90);
\draw[->][green] (0,-4.5-2)--(20,13.50);
\draw[->][green] (0,-5.0-2)--(20,13.00);

\draw [->](1,1-2)-- (20,18);
\node at (1,1-2) [fill=white, inner sep=1pt, shape=circle, draw] {};

\draw [->](2,0-2)-- (20,16);
\node at (2,0-2) [fill=black, inner sep=1pt, shape=circle, draw] {};

\draw[->][red] (2.5,0-2)--(20,15.50);
\draw[->][red] (2.5,-0.5-2)--(20,15.00);
\node at (2.5,0-2) [fill=red, inner sep=1pt, shape=circle, draw] {};
\node at (2.5,-0.5-2) [fill=red, inner sep=1pt, shape=circle, draw] {};

\draw [->](3,-1-2)-- (20,14.0);
\node at (3,-1-2) [fill=white, inner sep=1pt, shape=circle, draw] {};

\draw [->](4,-2-2)-- (20,12.0);
\node at (4,-2-2) [fill=white, inner sep=1pt, shape=diamond, draw] {};

\draw[->][red] (4.5,-2.0-2)--(20,11.5);
\draw[->][red] (4.5,-2.5-2)--(20,11.0);

\draw[->][blue] (5.5,-2.0-2)--(20,10.50);
\draw[->][blue] (5.5,-2.60-2)--(20,9.90);
\draw[->][blue] (5.5,-3.0-2)--(20,9.50);
\draw[->][blue] (5.5,-3.5-2)--(20,9.00);

\node at (4.5,-2.0-2) [fill=red, inner sep=1pt, shape=circle, draw] {};
\node at (4.5,-2.5-2) [fill=red, inner sep=1pt, shape=circle, draw] {};
\node at (5.5,-2.0-2) [fill=red, inner sep=1pt, shape=circle, draw] {};
\node at (5.5,-2.5-2) [fill=red, inner sep=1pt, shape=circle, draw] {};
\node at (5.5,-3.0-2) [fill=red, inner sep=1pt, shape=circle, draw] {};
\node at (5.5,-3.5-2) [fill=red, inner sep=1pt, shape=circle, draw] {};

\draw [->](5,-3-2)-- (20,10.0);
\node at (5,-3-2) [fill=white, inner sep=1pt, shape=circle, draw] {};

\draw [->](6,-4-2)-- (20,8.0);
\node at (6,-4-2) [fill=black, inner sep=1pt, shape=circle, draw] {};

\draw[->][red] (6.5,-4.0-2)--(20,7.5);
\draw[->][red] (6.5,-4.5-2)--(20,7.0);
\node at (6.5,-4.0-2) [fill=red, inner sep=1pt, shape=circle, draw] {};
\node at (6.5,-4.5-2) [fill=red, inner sep=1pt, shape=circle, draw] {};

\draw [->](7,-5-2)-- (20,6.0);
\node at (7,-5-2) [fill=white, inner sep=1pt, shape=circle, draw] {};

%%(BB2)
\foreach \y in {1,2,3,4,5,6,7,8,9,10,11,12,13,14,15,16,17,18,19,20,21,22,23,24}
\node at (8,-8-\y/2) [fill=red, inner sep=1pt, shape=circle, draw] {};

\draw [->] (8,-8)-- (20,20-16);
\node at (8,-8)  [shape = rectangle, draw]{};
\draw (8,-8) node[anchor=east]{$U^1$};

\draw[->][red] (8,-8-0.5)--(20,19.5-16);
\draw[->][red] (8,-8-1.0)--(20,19.0-16);

\draw[->][blue] (0+8,0.5-2-8)--(20,18.5-16);
\draw[->][blue] (0+8,-0.10-2-8)--(20,17.90-16);
\draw[->][blue] (0+8,-0.5-2-8)--(20,17.50-16);
\draw[->][blue] (0+8,-1.0-2-8)--(20,17.00-16);

\draw[->][green] (0+8,-1.5-2-8)--(20,16.50-16);
\draw[->][green] (0+8,-2.10-2-8)--(20,15.90-16);
\draw[->][green] (0+8,-2.60-2-8)--(20,15.40-16);
\draw[->][green] (0+8,-3.10-2-8)--(20,14.90-16);
\draw[->][green] (0+8,-3.5-2-8)--(20,14.50-16);
\draw[->][green] (0+8,-4.10-2-8)--(20,13.90-16);
\draw[->][green] (0+8,-4.5-2-8)--(20,13.50-16);
\draw[->][green] (0+8,-5.0-2-8)--(20,13.00-16);

\draw [->](1+8,1-2-8)-- (20,18-16);
\node at (1+8,1-2-8) [fill=white, inner sep=1pt, shape=circle, draw] {};

\draw [->](2+8,0-2-8)-- (20,16-16);
\node at (2+8,0-2-8) [fill=black, inner sep=1pt, shape=circle, draw] {};

\draw[->][red] (2.5+8,0-2-8)--(20,15.50-16);
\draw[->][red] (2.5+8,-0.5-2-8)--(20,15.00-16);
\node at (2.5+8,0-2-8) [fill=red, inner sep=1pt, shape=circle, draw] {};
\node at (2.5+8,-0.5-2-8) [fill=red, inner sep=1pt, shape=circle, draw] {};

\draw [->](3+8,-1-2-8)-- (20,14.0-16);
\node at (3+8,-1-2-8) [fill=white, inner sep=1pt, shape=circle, draw] {};

\draw [->](4+8,-2-2-8)-- (20,12.0-16);
\node at (4+8,-2-2-8) [fill=white, inner sep=1pt, shape=diamond, draw] {};

\draw[->][red] (4.5+8,-2.0-2-8)--(20,11.5-16);
\draw[->][red] (4.5+8,-2.5-2-8)--(20,11.0-16);

\draw[->][blue] (5.5+8,-2.0-2-8)--(20,10.50-16);
\draw[->][blue] (5.5+8,-2.60-2-8)--(20,9.90-16);
\draw[->][blue] (5.5+8,-3.0-2-8)--(20,9.50-16);
\draw[->][blue] (5.5+8,-3.5-2-8)--(20,9.00-16);

\node at (4.5+8,-2.0-2-8) [fill=red, inner sep=1pt, shape=circle, draw] {};
\node at (4.5+8,-2.5-2-8) [fill=red, inner sep=1pt, shape=circle, draw] {};
\node at (5.5+8,-2.0-2-8) [fill=red, inner sep=1pt, shape=circle, draw] {};
\node at (5.5+8,-2.5-2-8) [fill=red, inner sep=1pt, shape=circle, draw] {};
\node at (5.5+8,-3.0-2-8) [fill=red, inner sep=1pt, shape=circle, draw] {};
\node at (5.5+8,-3.5-2-8) [fill=red, inner sep=1pt, shape=circle, draw] {};

\draw [->](5+8,-3-2-8)-- (20,10.0-16);
\node at (5+8,-3-2-8) [fill=white, inner sep=1pt, shape=circle, draw] {};

\draw [->](6+8,-4-2-8)-- (20,8.0-16);
\node at (6+8,-4-2-8) [fill=black, inner sep=1pt, shape=circle, draw] {};

\draw[->][red] (6.5+8,-4.0-2-8)--(20,7.5-16);
\draw[->][red] (6.5+8,-4.5-2-8)--(20,7.0-16);
\node at (6.5+8,-4.0-2-8) [fill=red, inner sep=1pt, shape=circle, draw] {};
\node at (6.5+8,-4.5-2-8) [fill=red, inner sep=1pt, shape=circle, draw] {};

\draw [->](7+8,-5-2-8)-- (20,6.0-16);
\node at (7+8,-5-2-8) [fill=white, inner sep=1pt, shape=circle, draw] {};
%%(BB3)
\foreach \y in {1,2,3,4,5,6,7,8}
\node at (8+8,-8-8-\y/2) [fill=red, inner sep=1pt, shape=circle, draw] {};

\draw [->] (8+8,-8-8)-- (20,20-16-16);
\node at (8+8,-8-8)  [shape = rectangle, draw]{};
\draw (8+8,-8-8) node[anchor=east]{$U^2$};

\draw[->][red] (8+8,-8-8-0.5)--(20,19.5-16-16);
\draw[->][red] (8+8,-8-8-1.0)--(20,19.0-16-16);

\draw[->][blue] (0+8+8,0.5-2-8-8)--(20,18.5-16-16);
\draw[->][blue] (0+8+8,-0.10-2-8-8)--(20,17.90-16-16);
\draw[->][blue] (0+8+8,-0.5-2-8-8)--(20,17.50-16-16);
\draw[->][blue] (0+8+8,-1.0-2-8-8)--(20,17.00-16-16);

\draw[->][green] (0+8+8,-1.5-2-8-8)--(20,16.50-16-16);
\draw[->][green] (0+8+8,-2.10-2-8-8)--(20,15.90-16-16);
\draw[->][green] (0+8+8+0.5,-2.10-2-8-8)--(20,15.40-16-16);
\draw[->][green] (0+8+8+1.0,-2.10-2-8-8)--(20,14.90-16-16);
\draw[->][green] (0+8+8+1.5,-2.10-2-8-8)--(20,14.50-16-16);
\draw[->][green] (0+8+8+2.0,-2.10-2-8-8)--(20,13.90-16-16);
\draw[->][green] (0+8+8+2.5,-2.10-2-8-8)--(20,13.50-16-16);
\draw[->][green] (0+8+8+3.0,-2.10-2-8-8)--(20,13.00-16-16);

\draw [->](1+8+8,1-2-8-8)-- (20,18-16-16);
\node at (1+8+8,1-2-8-8) [fill=white, inner sep=1pt, shape=circle, draw] {};

\draw [->](2+8+8,0-2-8-8)-- (20,16-16-16);
\node at (2+8+8,0-2-8-8) [fill=black, inner sep=1pt, shape=circle, draw] {};

\draw[->][red] (2.5+8+8,0-2-8-8)--(20,15.50-16-16);
\draw[->][red] (2.5+8+8,-0.5-2-8-8)--(20,15.00-16-16);
\node at (2.5+8+8,0-2-8-8) [fill=red, inner sep=1pt, shape=circle, draw] {};
\node at (2.5+8+8,-0.5-2-8-8) [fill=red, inner sep=1pt, shape=circle, draw] {};

\draw [->](3+8+8,-1-2-8-8)-- (20,14.0-16-16);
\node at (3+8+8,-1-2-8-8) [fill=white, inner sep=1pt, shape=circle, draw] {};

\node at (4+8+8,-2-2-8-8) [fill=white, inner sep=1pt, shape=diamond, draw] {};
\foreach \y in {1,2,3,4,5,6,7,8,9,10,11,12,13,14,15,16,17,18,19,20,21,22,23,24}
\node at (-4.5-4,10-2+\y/2) [fill=red, inner sep=1pt, shape=circle, draw] {};

\draw [->] (3-4,3-2)-- (18,20);
\node at (3-4,3-2) [fill=white, inner sep=1pt, shape=circle, draw] {};

\draw[->][red] (2.5-4,3.5-2)--(17.0,20.0);
\draw[->][red] (2.5-4,4.0-2)--(16.5,20.0);
\node at (2.5-4,3.5-2) [fill=red, inner sep=1pt, shape=circle, draw] {};
\node at (2.5-4,4.0-2) [fill=red, inner sep=1pt, shape=circle, draw] {};

\draw [->] (2.0-4,4.0-2)-- (16.0,20);
\node at (2.0-4,4.0-2) [fill=black, inner sep=1pt, shape=circle, draw] {};

\draw[->][blue] (1.5-4,4.5-2)--(15.0,20);
\draw[->][blue] (1.5-4,5.0-2)--(14.5,20);
\draw[->][blue] (1.5-4,5.40-2)--(14.0,19.90);
\draw[->][blue] (1.5-4,6.0-2)--(13.5,20);
\node at (1.5-4,4.5-2) [fill=red, inner sep=1pt, shape=circle, draw] {};
\node at (1.5-4,5.0-2) [fill=red, inner sep=1pt, shape=circle, draw] {};
\node at (1.5-4,5.5-2) [fill=red, inner sep=1pt, shape=circle, draw] {};
\node at (1.5-4,6.0-2) [fill=red, inner sep=1pt, shape=circle, draw] {};

\draw [->] (1.0-4,5.0-2)-- (14.0,20);
\node at (1.0-4,5.0-2) [fill=white, inner sep=1pt, shape=circle, draw] {};

\draw [->] (0.0-4,6.0-2)-- (12.0,20);
\node at (0.0-4,6.0-2) [fill=white, inner sep=1pt, shape=diamond, draw] {};

\draw[->][red] (0.5-4,6.0-2)--(12.5,20);
\draw[->][red] (0.5-4,5.5-2)--(13.0,20);
\node at (0.5-4,6.0-2) [fill=red, inner sep=1pt, shape=circle, draw] {};
\node at (0.5-4,5.5-2) [fill=red, inner sep=1pt, shape=circle, draw] {};

\draw[->][green] (-0.5-4,6.5-2)--(11.0,20);
\draw[->][green] (-0.5-4,7.0-2)--(10.5,20);
\draw[->][green] (-0.5-4,7.40-2)--(10.0,19.90);
\draw[->][green] (-0.5-4,8.00-2)--(9.5,20);
\draw[->][green] (-0.5-4,8.40-2)--(9.0,19.90);
\draw[->][green] (-0.5-4,8.90-2)--(8.5,19.90);
\draw[->][green] (-0.5-4,9.40-2)--(8.0,19.90);
\draw[->][green] (-0.5-4,10.0-2)--(7.5,20);
\node at (-0.5-4,6.5-2) [fill=red, inner sep=1pt, shape=circle, draw] {};
\node at (-0.5-4,7.0-2) [fill=red, inner sep=1pt, shape=circle, draw] {};
\node at (-0.5-4,7.5-2) [fill=red, inner sep=1pt, shape=circle, draw] {};
\node at (-0.5-4,8.0-2) [fill=red, inner sep=1pt, shape=circle, draw] {};
\node at (-0.5-4,8.5-2) [fill=red, inner sep=1pt, shape=circle, draw] {};
\node at (-0.5-4,9.0-2) [fill=red, inner sep=1pt, shape=circle, draw] {};
\node at (-0.5-4,9.5-2) [fill=red, inner sep=1pt, shape=circle, draw] {};
\node at (-0.5-4,10.0-2) [fill=red, inner sep=1pt, shape=circle, draw] {};

\draw [->] (-1.0-4,7.0-2)-- (10.0,20);
\node at (-1.0-4,7.0-2) [fill=white, inner sep=1pt, shape=circle, draw] {};

\draw [->] (-2.0-4,8.0-2)-- (8.0,20);
\node at (-2.0-4,8.0-2) [fill=black, inner sep=1pt, shape=circle, draw] {};

\draw [->] (-3.0-4,9.0-2)-- (6.0,20);
\node at (-3.0-4,9.0-2) [fill=white, inner sep=1pt, shape=circle, draw] {};

\draw[->][red] (-1.5-4,8.0-2)--(8.5,20);
\draw[->][red] (-1.5-4,7.5-2)--(9.0,20);
\node at (-1.5-4,8.0-2) [fill=red, inner sep=1pt, shape=circle, draw] {};
\node at (-1.5-4,7.5-2) [fill=red, inner sep=1pt, shape=circle, draw] {};

\draw[->][blue] (-2.5-4,8.5-2)--(7.0,20);
\draw[->][blue] (-2.5-4,9.0-2)--(6.5,20);
\draw[->][blue] (-2.5-4,9.40-2)--(6.0,19.90);
\draw[->][blue] (-2.5-4,10.0-2)--(5.5,20);
\node at (-2.5-4,8.5-2) [fill=red, inner sep=1pt, shape=circle, draw] {};
\node at (-2.5-4,9.0-2) [fill=red, inner sep=1pt, shape=circle, draw] {};
\node at (-2.5-4,9.5-2) [fill=red, inner sep=1pt, shape=circle, draw] {};
\node at (-2.5-4,10.0-2) [fill=red, inner sep=1pt, shape=circle, draw] {};

\draw [->] (-4.0-4,10.0-2)-- (4.0,20);
\node at (-4.0-4,10.0-2) [inner sep=3pt, shape=circle, draw] {};
\draw (-4.0-4-0.25,10.0-2) node[anchor=east]{$2 \cdot U^{-1}$};

\draw[->][red] (-3.5-4,10.0-2)--(4.5,20);
\draw[->][red] (-3.5-4,9.5-2)--(5.0,20);
\node at (-3.5-4,10.0-2) [fill=red, inner sep=1pt, shape=circle, draw] {};
\node at (-3.5-4,9.5-2) [fill=red, inner sep=1pt, shape=circle, draw] {};

\foreach \y in {1,2,3,4,5,6,7,8}
\node at (-4.5-4-8,10-2+8+\y/2) [fill=red, inner sep=1pt, shape=circle, draw] {};

\draw [->] (3-4-8,3-2+8)-- (18-16,20);
\node at (3-4-8,3-2+8) [fill=white, inner sep=1pt, shape=circle, draw] {};

\draw[->][red] (2.5-4-8,3.5-2+8)--(17.0-16,20.0);
\draw[->][red] (2.5-4-8,4.0-2+8)--(16.5-16,20.0);
\node at (2.5-4-8,3.5-2+8) [fill=red, inner sep=1pt, shape=circle, draw] {};
\node at (2.5-4-8,4.0-2+8) [fill=red, inner sep=1pt, shape=circle, draw] {};

\draw [->] (2.0-4-8,4.0-2+8)-- (16.0-16,20);
\node at (2.0-4-8,4.0-2+8) [fill=black, inner sep=1pt, shape=circle, draw] {};

\draw[->][blue] (1.5-4-8,4.5-2+8)--(15.0-16,20);
\draw[->][blue] (1.5-4-8,5.0-2+8)--(14.5-16,20);
\draw[->][blue] (1.5-4-8,5.40-2+8)--(14.0-16,19.90);
\draw[->][blue] (1.5-4-8,6.0-2+8)--(13.5-16,20);
\node at (1.5-4-8,4.5-2+8) [fill=red, inner sep=1pt, shape=circle, draw] {};
\node at (1.5-4-8,5.0-2+8) [fill=red, inner sep=1pt, shape=circle, draw] {};
\node at (1.5-4-8,5.5-2+8) [fill=red, inner sep=1pt, shape=circle, draw] {};
\node at (1.5-4-8,6.0-2+8) [fill=red, inner sep=1pt, shape=circle, draw] {};

\draw [->] (1.0-4-8,5.0-2+8)-- (14.0-16,20);
\node at (1.0-4-8,5.0-2+8) [fill=white, inner sep=1pt, shape=circle, draw] {};

\draw [->] (0.0-4-8,6.0-2+8)-- (12.0-16,20);
\node at (0.0-4-8,6.0-2+8) [fill=white, inner sep=1pt, shape=diamond, draw] {};

\draw[->][red] (0.5-4-8,6.0-2+8)--(12.50-16,20);
\draw[->][red] (0.5-4-8,5.5-2+8)--(13.00-16,20);
\node at (0.5-4-8,6.0-2+8) [fill=red, inner sep=1pt, shape=circle, draw] {};
\node at (0.5-4-8,5.5-2+8) [fill=red, inner sep=1pt, shape=circle, draw] {};

\draw[->][green] (-0.5-4-8,6.5-2+8)--(11.0-16,20);
\draw[->][green] (-0.5-4-8,7.0-2+8)--(10.5-16,20);
\draw[->][green] (-0.5-4-8,7.40-2+8)--(10.0-16,19.90);
\draw[->][green] (-0.5-4-8,8.00-2+8)--(9.5-16,20.00);
\draw[->][green] (-0.5-4-8,8.40-2+8)--(9.0-16,19.90);
\draw[->][green] (-0.5-4-8,8.90-2+8)--(8.5-16,19.90);
\draw[->][green] (-0.5-4-8,9.40-2+8)--(8.0-16,19.90);
\draw[->][green] (-0.5-4-8,10.0-2+8)--(7.5-16,20);
\node at (-0.5-4-8,6.5-2+8) [fill=red, inner sep=1pt, shape=circle, draw] {};
\node at (-0.5-4-8,7.0-2+8) [fill=red, inner sep=1pt, shape=circle, draw] {};
\node at (-0.5-4-8,7.5-2+8) [fill=red, inner sep=1pt, shape=circle, draw] {};
\node at (-0.5-4-8,8.0-2+8) [fill=red, inner sep=1pt, shape=circle, draw] {};
\node at (-0.5-4-8,8.5-2+8) [fill=red, inner sep=1pt, shape=circle, draw] {};
\node at (-0.5-4-8,9.0-2+8) [fill=red, inner sep=1pt, shape=circle, draw] {};
\node at (-0.5-4-8,9.5-2+8) [fill=red, inner sep=1pt, shape=circle, draw] {};
\node at (-0.5-4-8,10.0-2+8) [fill=red, inner sep=1pt, shape=circle, draw] {};

\draw [->] (-1.0-4-8,7.0-2+8)-- (10.0-16,20);
\node at (-1.0-4-8,7.0-2+8) [fill=white, inner sep=1pt, shape=circle, draw] {};

\draw [->] (-2.0-4-8,8.0-2+8)-- (8.0-16,20);
\node at (-2.0-4-8,8.0-2+8) [fill=black, inner sep=1pt, shape=circle, draw] {};

\draw [->] (-3.0-4-8,9.0-2+8)-- (6.0-16,20);
\node at (-3.0-4-8,9.0-2+8) [fill=white, inner sep=1pt, shape=circle, draw] {};

\draw[->][red] (-1.5-4-8,8.0-2+8)--(8.5-16,20);
\draw[->][red] (-1.5-4-8,7.5-2+8)--(9.0-16,20);
\node at (-1.5-4-8,8.0-2+8) [fill=red, inner sep=1pt, shape=circle, draw] {};
\node at (-1.5-4-8,7.5-2+8) [fill=red, inner sep=1pt, shape=circle, draw] {};

\draw[->][blue] (-2.5-4-8,8.5-2+8)--(7.0-16,20);
\draw[->][blue] (-2.5-4-8,9.0-2+8)--(6.5-16,20);
\draw[->][blue] (-2.5-4-8,9.40-2+8)--(6.0-16,19.90);
\draw[->][blue] (-2.5-4-8,10.0-2+8)--(5.5-16,20);
\node at (-2.5-4-8,8.5-2+8) [fill=red, inner sep=1pt, shape=circle, draw] {};
\node at (-2.5-4-8,9.0-2+8) [fill=red, inner sep=1pt, shape=circle, draw] {};
\node at (-2.5-4-8,9.5-2+8) [fill=red, inner sep=1pt, shape=circle, draw] {};
\node at (-2.5-4-8,10.0-2+8) [fill=red, inner sep=1pt, shape=circle, draw] {};

\draw [->] (-4.0-4-8,10.0-2+8)-- (4.0-16,20);
\node at (-4.0-4-8,10.0-2+8) [inner sep=3pt, shape=circle, draw] {};
\draw (-4.0-4-8-0.25,10.0-2+8) node[anchor=east]{$2 \cdot U^{-2}$};

\draw[->][red] (-3.5-4-8,10.0-2+8)--(4.5-16,20);
\draw[->][red] (-3.5-4-8,9.5-2+8)--(5.0-16,20);
\node at (-3.5-4-8,10.0-2+8) [fill=red, inner sep=1pt, shape=circle, draw] {};
\node at (-3.5-4-8,9.5-2+8) [fill=red, inner sep=1pt, shape=circle, draw] {};

\draw [->] (3-4-8-8,3-2+8+8)-- (18-16-16,20);
\node at (3-4-8-8,3-2+8+8) [fill=white, inner sep=1pt, shape=circle, draw] {};

\draw[->][red] (2.5-4-8-8,3.5-2+8+8)--(17.0-16-16,20.0);
\draw[->][red] (2.5-4-8-8,4.0-2+8+8)--(16.5-16-16,20.0);
\node at (2.5-4-8-8,3.5-2+8+8) [fill=red, inner sep=1pt, shape=circle, draw] {};
\node at (2.5-4-8-8,4.0-2+8+8) [fill=red, inner sep=1pt, shape=circle, draw] {};

\draw [->] (2.0-4-8-8,4.0-2+8+8)-- (16.0-16-16,20);
\node at (2.0-4-8-8,4.0-2+8+8) [fill=black, inner sep=1pt, shape=circle, draw] {};

\draw[->][blue] (1.5-4-8-8,4.5-2+8+8)--(15.0-16-16,20);
\draw[->][blue] (1.5-4-8-8,5.0-2+8+8)--(14.5-16-16,20);
\draw[->][blue] (1.5-4-8-8,5.40-2+8+8)--(14.0-16-16,19.90);
\node at (1.5-4-8-8,4.5-2+8+8) [fill=red, inner sep=1pt, shape=circle, draw] {};
\node at (1.5-4-8-8,5.0-2+8+8) [fill=red, inner sep=1pt, shape=circle, draw] {};
\node at (1.5-4-8-8,5.5-2+8+8) [fill=red, inner sep=1pt, shape=circle, draw] {};
\node at (1.5-4-8-8,6.0-2+8+8) [fill=red, inner sep=1pt, shape=circle, draw] {};

\draw [->] (1.0-4-8-8,5.0-2+8+8)-- (14.0-16-16,20);
\node at (1.0-4-8-8,5.0-2+8+8) [fill=white, inner sep=1pt, shape=circle, draw] {};

\node at (0.0-4-8-8,6.0-2+8+8) [fill=white, inner sep=1pt, shape=diamond, draw] {};

\draw[->][red] (0.5-4-8-8,5.5-2+8+8)--(13.00-16-16,20);
\node at (0.5-4-8-8,6.0-2+8+8) [fill=red, inner sep=1pt, shape=circle, draw] {};
\node at (0.5-4-8-8,5.5-2+8+8) [fill=red, inner sep=1pt, shape=circle, draw] {};

\draw[->][green] (-0.5-4-8,6.5-2+8)--(11.0-16,20);
\draw[->][green] (-0.5-4-8,7.0-2+8)--(10.5-16,20);
\draw[->][green] (-0.5-4-8,7.40-2+8)--(10.0-16,19.90);
\draw[->][green] (-0.5-4-8,8.00-2+8)--(9.5-16,20.00);
\draw[->][green] (-0.5-4-8,8.40-2+8)--(9.0-16,19.90);
\draw[->][green] (-0.5-4-8,8.90-2+8)--(8.5-16,19.90);
\draw[->][green] (-0.5-4-8,9.40-2+8)--(8.0-16,19.90);
\draw[->][green] (-0.5-4-8,10.0-2+8)--(7.5-16,20);
\node at (-0.5-4-8,6.5-2+8) [fill=red, inner sep=1pt, shape=circle, draw] {};
\node at (-0.5-4-8,7.0-2+8) [fill=red, inner sep=1pt, shape=circle, draw] {};
\node at (-0.5-4-8,7.5-2+8) [fill=red, inner sep=1pt, shape=circle, draw] {};
\node at (-0.5-4-8,8.0-2+8) [fill=red, inner sep=1pt, shape=circle, draw] {};
\node at (-0.5-4-8,8.5-2+8) [fill=red, inner sep=1pt, shape=circle, draw] {};
\node at (-0.5-4-8,9.0-2+8) [fill=red, inner sep=1pt, shape=circle, draw] {};
\node at (-0.5-4-8,9.5-2+8) [fill=red, inner sep=1pt, shape=circle, draw] {};
\node at (-0.5-4-8,10.0-2+8) [fill=red, inner sep=1pt, shape=circle, draw] {};

\draw[->][red] (-1.5-4-8,8.0-2+8)--(8.5-16,20);
\draw[->][red] (-1.5-4-8,7.5-2+8)--(9.0-16,20);
\node at (-1.5-4-8,8.0-2+8) [fill=red, inner sep=1pt, shape=circle, draw] {};
\node at (-1.5-4-8,7.5-2+8) [fill=red, inner sep=1pt, shape=circle, draw] {};

\draw[->][blue] (-2.5-4-8,8.5-2+8)--(7.0-16,20);
\draw[->][blue] (-2.5-4-8,9.0-2+8)--(6.5-16,20);
\draw[->][blue] (-2.5-4-8,9.40-2+8)--(6.0-16,19.90);
\draw[->][blue] (-2.5-4-8,10.0-2+8)--(5.5-16,20);
\node at (-2.5-4-8,8.5-2+8) [fill=red, inner sep=1pt, shape=circle, draw] {};
\node at (-2.5-4-8,9.0-2+8) [fill=red, inner sep=1pt, shape=circle, draw] {};
\node at (-2.5-4-8,9.5-2+8) [fill=red, inner sep=1pt, shape=circle, draw] {};
\node at (-2.5-4-8,10.0-2+8) [fill=red, inner sep=1pt, shape=circle, draw] {};

\draw [->] (-4.0-4-8,10.0-2+8)-- (4.0-16,20);
\node at (-4.0-4-8,10.0-2+8) [inner sep=3pt, shape=circle, draw] {};
\draw (-4.0-4-8-0.25,10.0-2+8) node[anchor=east]{$2 \cdot U^{-2}$};

\draw[->][red] (-3.5-4-8,10.0-2+8)--(4.5-16,20);
\draw[->][red] (-3.5-4-8,9.5-2+8)--(5.0-16,20);
\node at (-3.5-4-8,10.0-2+8) [fill=red, inner sep=1pt, shape=circle, draw] {};
\node at (-3.5-4-8,9.5-2+8) [fill=red, inner sep=1pt, shape=circle, draw] {};

\end{tikzpicture}$$
\caption{The homotopy structure of $BP\mathbb R\langle 3 \rangle^{Q}_{\rost}$} \label{fig:GL3}
\end{figure}

\bigskip

\section{Notation}
We adopt the notation from \cite{BPRn}.  We write $\Z$ for the 2-local integers. The element $u$ has degree
$2(1-\sigma)$, with $U=u^8$ and  $\vb_i$ is an element of degree
$(2^i-1)\rho$. We will be working with the ring
$$P=\Z[\vb_1, \vb_2, \ldots , \vb_n]$$
and its ideal
$$J=(\vb_1, \ldots, \vb_n). $$
We also need to consider its quotients
$$P_i=\Z [\vb_{i+1}, \ldots , \vb_n]$$
and similarly for the mod 2 reductions
$$\Pb =\Ftwo [\vb_1, \vb_2, \ldots , \vb_n]$$
and its quotients
$$\Pb_i=\Ftwo  [\vb_{i+1}, \ldots , \vb_n]. $$
The $J$-local cohomology of all these $P$-modules is easy to
determine, since they are all $J$-Cohen-Macaulay (so the  local
cohomology is all in the top degree) and $J$-Gorenstein so that
the top local cohomology is the dual (in the appropriate sense) of the
original ring:
$$H^*_J(P_i)=H^{n-i}_J(P_i)=
\Sigma^{|\vb_{i+1}|+\cdots +|\vb_n|} P_i^*, \mbox{ where } P_i^*=\Hom_{\Z} (P_i, \Z) $$
and
$$H^*_J(\Pb_i)=H^{n-i}_J(\Pb_i)=
\Sigma^{|\vb_{i+1}|+\cdots +|\vb_n|}\Pb_i^{\vee} \mbox{
  where } \Pb_i^{\vee}=\Hom_{\Ftwo}(\Pb_i, \Ftwo). $$
This lets one deduce the local cohomology of the kernel of a
surjective map from one of these to another.

We will be exclusively restricting attention to $n=3$ so we make explicit
\begin{itemize}
\item $P = BP\mathbb{R}\langle 3 \rangle_{* \rho}^{C_2} = \Z  [\bar v_1, \bar v_2, \bar v_3]$,
\item $\Pb_0 = \mathbb{F}_2 [\bar v_1, \bar v_2, \bar v_3]$,
\item $\Pb_1 = \mathbb{F}_2 [\bar v_2, \bar v_3]$,
\item $\Pb_2 = \mathbb{F}_2 [\bar v_3]$ and
\item $\Pb_3 = \mathbb{F}_2$.
\end{itemize}

Referring back to Figure \ref{fig:GL3},
\begin{itemize}
\item the black diagonal lines which start in a square indicate the copies of $P$,
\item the black diagonal lines which start in a small circle indicate the copies of $(2)P$,
\item the black diagonal lines which start in a dot indicate the copies of $(2,\bar v_1)P$,
\item the black diagonal lines which start in a diamond indicate the copies of $(2,\bar v_1,\bar v_2)P$,
\item the black diagonal lines which start in a big circle indicate the copies of $(2,\bar v_1,\bar v_2,\bar v_3)P$,
\item the red diagonal lines indicate the copies of $\Pb_0$,
\item the blue diagonal lines indicate the copies of $\Pb_1$,
\item the green diagonal lines indicate the copies of $\Pb_2$, and
\item the red dots are the copies of $\mathbb{F}_2 = \Pb_3$.
\end{itemize}

\bigskip

\section{Local cohomology and the spectral sequence}

We next wish to calculate the $J$-local cohomology. Since $J$ is
generated by elements in degrees which are multiples of $\rho$, the
local cohomology will respect the decomposition of a module into its
diagonal pieces (i.e., into the $\delta$-diagonals, consisting of elements
of degrees  $\delta+*\rho$).

We have the following tables of the basic block (Table
\ref{tab:table1}) and the local cohomology (Tables \ref{tab:table2}
and \ref{tab:table3}) based on $BP \mathbb R \langle 3
\rangle^{Q}_{\rost}$ as follows. The number $\delta$ refers to the diagonal
(i.e., $\delta$ is the point on the $\Z$ axis where $\delta+*\rho$
intersects the $y=0$ line). The column refers to the  power of $u$ 
for that module (this power of $u$ is usually the bottom of the diagonal, but
not always). The modules themselves are self-explanatory.

\newpage

\renewcommand{\tabcolsep}{15.00pt}
\renewcommand{\arraystretch}{1.00}
\newcolumntype{a}{>{\columncolor{Gray}}c}
%\newcolumntype{b}{>{\columncolor{white}}c}
\begin{table}[h!]
\begin{center}
 \caption{Basic block of $BP \mathbb R\langle 3\rangle^{Q}_{\rost}$}
    \label{tab:table1}
    \begin{tabular}{ ||!{\mystrut}a | p{0.35cm} | p{0.50cm} | p{1.0cm} |p{0.50cm} | p{1.6cm}|p{0.50cm} |p{1.0cm} |p{0.50cm} || }
    \hline
      $\delta$ &\cellcolor{yellow!55}$1$ &\cellcolor{yellow!55}$u$ &\cellcolor{yellow!55}$u^2$ &\cellcolor{yellow!55}$u^3$ &\cellcolor{yellow!55}$u^4$
      &\cellcolor{yellow!55}$u^5$ &\cellcolor{yellow!55}$u^6$ &\cellcolor{yellow!55}$u^7$  \\ [1.5ex]
    \hline \hline
      $0$ &$P$         &       &            &       & & & &\\    \hline
      $1$ &$\Pb_0$  &       &            &       & & & &\\    \hline
      $2$ &$\Pb_0$  &       &            &       & & & &\\    \hline
      $3$ &$\Pb_1$  &       &            &       & & & &\\    \hline
      $4$ &$\Pb_1$  &$(2)P$ &            &       & & & &\\    \hline
      $5$ &$\Pb_1$  &       &            &       & & & &\\    \hline
      $6$ &$\Pb_1$  &       &            &       & & & &\\    \hline
      $7$ &$\Pb_2$  &       &            &       & & & &\\    \hline
      $8$ &$\Pb_2$  &       &$(2,\bar v_1)P$  &       & & & &\\    \hline
      $9$ &$\Pb_2$  &       &$(\bar v_1)\Pb_0$  &       & & & &\\    \hline
      $10$ &$\Pb_2$ &       &$(\bar v_1)\Pb_0$  &       & & & &\\    \hline
      $11$ &$\Pb_2$ &       &            &       & & & &\\    \hline
      $12$ &$\Pb_2$ &       &            &$(2)P$ &      & & &\\ \hline
      $13$ &$\Pb_2$ &       &            &       &      & & &\\ \hline
      $14$ &$\Pb_2$ &       &            &       &      & & &\\ \hline
      $15$ &$\Pb_3$ &       &            &       &      & & &\\ \hline
      $16$ &$\Pb_3$ &       &            &       &$(2,\bar v_1,\bar v_2)P$      & & &\\ \hline
      $17$ &$\Pb_3$ &       &            &       &$(\bar v_1,\bar v_2)\Pb_0$ & & &\\ \hline
      $18$ &$\Pb_3$ &       &            &       &$(\bar v_1,\bar v_2)\Pb_0$ & & &\\ \hline
      $19$ &$\Pb_3$ &       &            &       &$(\bar v_2)\Pb_1$ & & &\\ \hline
      $20$ &$\Pb_3$ &       &            &       &$(\bar v_2)\Pb_1$ &$(2)P$  & &\\ \hline
      $21$ &$\Pb_3$ &       &            &       &$(\bar v_2)\Pb_1$ &       & &\\ \hline
      $22$ &$\Pb_3$ &       &            &       &$(\bar v_2)\Pb_1$ &       & &\\ \hline
      $23$ &$\Pb_3$ &       &            &       &                    &       & &\\ \hline
      $24$ &$\Pb_3$ &       &            &       &                    & &$(2,\bar v_1)P$ &\\ \hline
      $25$ &$\Pb_3$ &       &            &       &                    & &$(\bar v_1)\Pb_0$ &\\ \hline
      $26$ &$\Pb_3$ &       &            &       &                    & &$(\bar v_1)\Pb_0$ &\\ \hline
      $27$ &$\Pb_3$ &       &            &       &                    & &           &\\ \hline
      $28$ &$\Pb_3$ &       &            &       &                    & &           &$(2)P$\\ \hline
      $29$ &$\Pb_3$ &       &            &       &                    & &           &\\ \hline
      $30$ &$\Pb_3$ &       &            &       &                    & &           &\\ \hline
      $\delta \geq 31$ &$\Pb_3$ &       &            &       &                    & &           &\\
    \hline
    \end{tabular}
\end{center}
\end{table}
\vspace{0.2cm}

The local cohomology of each module is placed in the same $u$-power
column, but because the local cohomological degree shifts left, this
is in a lower $\delta$-row than the original module. Within that diagonal,
there is a shift from $u^i$ itself, indicated by the number in
parentheses. For example,  consider the entry $P$ in row 0 and column
$u^0$ of Table 1.  We calculate $H^*_J(P)=H^3_J(P)=\Sigma^{-11\rho}P^*$
($11\rho=|\vb_1|+|\vb_2|+|\vb_3|$), and so
we see an entry $P^*(-14\rho)$ (where $14=11+3$) in row $-3=0-3$, column
$u^0$ of Table 2. 

%\newpage

\renewcommand{\tabcolsep}{18.00pt}
\renewcommand{\arraystretch}{0.50}
\newcolumntype{a}{>{\columncolor{Gray}}c}
\begin{table}[h!]
\begin{center}
 \caption{Local cohomology of the basic block for $1, u, u^2$ and $u^3$}
    \label{tab:table2}
    \begin{tabular}{ ||!{\mystrut}a | p{1.65cm} | p{1.65cm} | p{4.20cm} |p{1.65cm} ||  }
    \hline
      $\delta$ &\cellcolor{yellow!55}$1$ &\cellcolor{yellow!55}$u$ &\cellcolor{yellow!55}$u^2$ &\cellcolor{yellow!55}$u^3$   \\ [1.5ex]
    \hline  \hline
      $-3$ &\textcolor{red}{$P^*(-14\rho)$}            &       &            &       \\    \hline
      $-2$ &\textcolor{red}{$\Pb_0^\vee(-14\rho)$}  &       &            &       \\    \hline
      $-1$ &\textcolor{red}{$\Pb_0^\vee(-14\rho)$}  &       &            &       \\    \hline
      $0$  &          &       &            &       \\    \hline
      $1$ &\textcolor{blue}{$\Pb_1^\vee(-12\rho)$}  &\textcolor{red}{$P^*(-14\rho)$} &            &      \\    \hline
      $2$ &\textcolor{blue}{$\Pb_1^\vee(-12\rho)$}  &       &            &      \\    \hline
      $3$ &\textcolor{blue}{$\Pb_1^\vee(-12\rho)$}  &       &            &      \\    \hline
      $4$ &\textcolor{blue}{$\Pb_1^\vee(-12\rho)$}  &       &            &       \\    \hline
      $5$ &  &       &\textcolor{red}{$P^*(-14\rho)$} $\oplus$ \textcolor{blue}{$\Pb_1^\vee(-13\rho)$}            &       \\    \hline
      $6$ &\textcolor{green}{$\Pb_2^\vee(-8\rho), d_2$}  &       &\textcolor{red}{$\Pb_0^\vee(-13\rho)$}            &       \\    \hline
      $7$ &\textcolor{green}{$\Pb_2^\vee(-8\rho), d_2$}  &       &\textcolor{red}{$\Pb_0^\vee(-13\rho)$}            &       \\    \hline
      $8$ &\textcolor{green}{$\Pb_2^\vee(-8\rho), d_2$}  &       &  &       \\    \hline
      $9$ &\textcolor{green}{$\Pb_2^\vee(-8\rho)$}  &       &  &\textcolor{red}{$P^*(-14\rho)$}       \\    \hline
      $10$ &\textcolor{green}{$\Pb_2^\vee(-8\rho)$} &       &  &       \\    \hline
      $11$ &\textcolor{green}{$\Pb_2^\vee(-8\rho)$} &       &            &       \\    \hline
      $12$ &\textcolor{green}{$\Pb_2^\vee(-8\rho)$} &       &            & \\ \hline
      $13$ &\textcolor{green}{$\Pb_2^\vee(-8\rho)$} &       &            &       \\ \hline
      $14$ &                                          &       &            &       \\ \hline
      $\delta \geq 15$ &$\Pb_3$ &       &            &      \\
      %$16$ &$\Pb_3$ &       &            &       \\ \hline
      %$17$ &$\Pb_3$ &       &            &      \\ \hline
      %$18$ &$\Pb_3$ &       &            &       \\ \hline
      %$19$ &$\Pb_3$ &       &            &      \\ \hline
      %$20$ &$\Pb_3$ &       &            &       \\ \hline
      %$21$ &$\Pb_3$ &       &            &       \\ \hline
      %$22$ &$\Pb_3$ &       &            &      \\ \hline
      %$23$ &$\Pb_3$ &       &            &      \\ \hline
      %$24$ &$\Pb_3$ &       &            &       \\ \hline
      %$25$ &$\Pb_3$ &       &            &       \\ \hline
      %$26$ &$\Pb_3$ &       &            &       \\ \hline
      %$27$ &$\Pb_3$ &       &            &       \\ \hline
      %$28$ &$\Pb_3$ &       &            &      \\ \hline
      %$29$ &$\Pb_3$ &       &            &       \\ \hline
      %$30$ &$\Pb_3$ &       &            &      \\ \hline
      %$31$ &$\Pb_3$ &       &            &      \\
    \hline
    \end{tabular}
\end{center}
\end{table}
\vspace{0.6cm}

%\newpage

\renewcommand{\tabcolsep}{13.48pt}
\renewcommand{\arraystretch}{0.5}
\newcolumntype{a}{>{\columncolor{Gray}}c}
\begin{table}[h!]
\begin{center}
 \caption{Local cohomology of the basic block for $u^4, u^5, u^6$ and $u^7$}
    \label{tab:table3}
    \begin{tabular}{ ||!{\mystrut}a | p{4.00cm} | p{1.65cm} | p{4.20cm} |p{1.65cm} ||  }
    \hline
      $\delta$ &\cellcolor{yellow!55}$u^4$ &\cellcolor{yellow!55}$u^5$ &\cellcolor{yellow!55}$u^6$ &\cellcolor{yellow!55}$u^7$   \\ [1.5ex]
    \hline \hline
      %$-3$ &            &       &            &       \\    \hline
      %$-2$ &  &       &            &       \\    \hline
      %$-1$ &  &       &            &       \\    \hline
      %$0$ &          &       &            &       \\    \hline
      %$1$ &  &  &            &      \\    \hline
      %$2$ &  &       &            &      \\    \hline
      %$3$ &  &       &            &      \\    \hline
      %$4$ &  &       &            &       \\    \hline
      %$5$ &  &       &            &       \\    \hline
      %$6$ &  &       &           &       \\    \hline
      %$7$ &  &       &            &       \\    \hline
      %$8$ &  &       &           &       \\    \hline
      %$9$ &  &       &           &       \\    \hline
      %$10$ &  &       &           &       \\    \hline
      %$11$ &  &       &           &       \\    \hline
      %$12$ &  &       &           &       \\    \hline
      $13$ &\textcolor{red}{$P^*(-14\rho)$}  &       &  &  \\    \hline
      $14$ &\textcolor{red}{$\Pb_0^\vee(-14\rho)$} $\oplus$ \textcolor{blue}{$\Pb_2^\vee(-9\rho)$} &       &  &       \\    \hline
      $15$ &\textcolor{red}{$\Pb_0^\vee(-14\rho)$} $\oplus$ \textcolor{blue}{$\Pb_2^\vee(-9\rho)$} &       &            &       \\    \hline
      $16$ &\hspace*{13ex}\textcolor{blue}{$\Pb_2^\vee(-9\rho)$} &       &            & \\ \hline
      $17$ &\textcolor{blue}{$\Pb_1^\vee(-10\rho)$} &\textcolor{red}{$P^*(-14\rho)$}       &            &       \\ \hline
      $18$ &\textcolor{blue}{$\Pb_1^\vee(-10\rho)$} &       &            &       \\ \hline
      $19$ &\textcolor{blue}{$\Pb_1^\vee(-10\rho)$} &       &            &      \\ \hline
      $20$ &\textcolor{blue}{$\Pb_1^\vee(-10\rho)$} &       &            &       \\ \hline
      $21$ & &       &\textcolor{red}{$P^*(-14\rho)$} $\oplus$ \textcolor{blue}{$\Pb_1^\vee(-13\rho)$}            &      \\ \hline
      $22$ & &       &\textcolor{red}{$\Pb_0^\vee(-13\rho)$}            &       \\ \hline
      $23$ & &       &\textcolor{red}{$\Pb_0^\vee(-13\rho)$}            &      \\ \hline
      $24$ & &       &            &       \\ \hline
      $25$ & &       &            &\textcolor{red}{$P^*(-14\rho)$}       \\ \hline
      $26$ & &       &            &      \\ \hline
      $27$ & &       &            &      \\ \hline
      $28$ & &       &            &       \\
      %$29$ & &       &            &       \\ \hline
      %$30$ & &       &            &       \\ \hline
      %$31$ & &       &            &       \\
    \hline
    \end{tabular}
\end{center}
\end{table}
\vspace{0.6cm}
Here we have coloured $H^1_J$-groups in green, $H^2_J$-groups in blue
and  $H^3_J$-groups in red, but note for example that the short exact sequence
$$0\lra (\vb_1,\vb_2)\lra \Pb \lra \Pb_2\lra 0$$ 
gives 
$$H^*_J((\vb_1, \vb_2))=H^3_J((\vb_1, \vb_2)) \oplus H^2_J((\vb_1,
\vb_2))\cong  H^3_J(\Pb) \oplus H^1_J(\Pb_2). $$
The colour coding makes the $H^1_J(\Pb_2)$ green, even though it
contributes to $H^2_J((\vb_1, \vb_2))$.

%\newpage

Figure \ref{fig:GL1} shows the basic block $BB$, the spectral sequence of the local cohomology and its Anderson dual. Here, we have mostly
omitted dots, circles and squares except at the ends of diagonals or where an additional generator is required. The vertical lines denote
multiplication by $a$. The yellow diamond does not denote a homotopy
class, but denotes the point we have to reflect (the non-torsion
classes) in to see Anderson duality. The torsion classes are shifted
by $-1$ after reflection; that is, shifted one step horizontally to
the left.

The symbol $GBB$ refers to the part of $\Gamma_JBP\R \langle 3 
\rangle$ coming from $BB$, and $SS(GBB)$ is the local cohomology in 
the $E^2$-term that calculates its homotopy. In other words, to go from 
$SS(GBB)$ to $GBB$ there may be differentials and additive 
extensions. The orange labels $a^i$ are the powers of $a$
that are in $H^0_J$, and otherwise the colours in $SS(GBB)$ correspond
to the colours of the diagonals in $BB$ that gave rise to them. 

Note that the Anderson dual of $SS(GBB)$ has been displayed on the (false)
assumption that there are no differentials and no additive extensions
in the  local cohomology. The point of this is that as displayed in
Figure \ref{fig:GL1}  we can see where classes must cancel to make the
Anderson dual a connective spectrum.

%\newpage

\begin{figure}
\centering
$$\begin{tikzpicture}[scale =0.705]
%\clip (-9, -8.6) rectangle (5, 8.5);
\draw[step=0.5, gray, very thin] (-12,-14) grid (12, 14);
\draw (5.0,-6.0 ) node[anchor=east, draw=orange]{\Large{BB}};
\draw (-8.0,-3.0 ) node[anchor=east, draw=orange]{\Large{SS(GBB)}};
\draw (-5.0,3.50 ) node[anchor=east, draw=orange]{\Large{$\mathbb{Z}^{\rm{SS(GBB)}}$}};

\foreach \y in {1,2,3,4,5,6,7,8,9,10,11,12,13,14}
\draw (0,2-\y/2) node[anchor=east] {\small{$a^{\y}$}};
\foreach \y in {15,16,17,18,19,20,21,22,23,24,25,26,27,28,29,30,31,32}
\draw (0,2-\y/2) node[anchor=west] {\color{orange}\small{$a^{\y}$}};

\foreach \y in {1,2,3,4,5,6,7,8,9,10,11,12,13,14,15,16,17,18,19,20,21,22,23,24,25,26,27,28,29,30,31,32}
\node at (0,2-\y/2) [fill=red, inner sep=1pt, shape=circle, draw] {};

\foreach \y in {0, 1,2,3,4,5,6,7,8,9,10,11,12,13,14,15,16,17,18,19,20,21,22,23}
\node at (\y/2+1/2,2+\y/2) [fill=red, inner sep=1pt, shape=circle, draw]{};

\foreach \y in {0, 1,2,3,4,5,6,7,8,9,10,11,12, 13,14,15,16,17,18,19,20,21,22,23}
\node at (\y/2+1/2,2+\y/2-1/2) [fill=red, inner sep=1pt, shape=circle, draw] {};

\draw [->] (0,2)-- (12,14);
\node at (0,2)  [shape = rectangle, draw]{};
\draw (0,2) node[anchor=east]{1};
\draw[red] (0,2)--(7,-5); \draw[red] (-1,3)--(-8,10);

\foreach \y in {1,2,3,4,5,6,7,8,9,10,11,12,13,14,15,16,17,18,19,20,21,22,23,24}
\draw (\y/2,2+\y/2) node[anchor=east] {\small{$\bar v_1^{\y}$}};

\foreach \y in {1,2,3,4,5,6,7,8,9,10,11,12,13,14,15,16,17,18,19,20,21,22,23,24}
\node at (\y/2,2+\y/2) [shape=rectangle, draw] {};

\draw [->] (-7,-3.5)-- (-12,-8.5);
\foreach \y in {-14,-15,-16,-17,-18,-19,-20,-21,-22,-23,-24}
\node at (\y/2,3.5+\y/2) [shape=rectangle, draw] {};

\draw [->] (-1,3)-- (10,14);
\node at (-1,3)  [shape = rectangle, draw]{};

\draw[->][red] (-1.5,3.5)--(9.0,14);
\draw[->][red] (-1.5,4.0)--(8.5,14);
\draw[red] (-1.5,3.5)--(-1.5,4.0);
\node at (-1.5,3.5) [fill=red, inner sep=1pt, shape=circle, draw] {};
\node at (-1.5,4.0) [fill=red, inner sep=1pt, shape=circle, draw] {};

\draw [->] (-2.0,4.0)-- (8.0,14);
\node at (-2.0,4.0)  [shape = circle, draw]{};

\draw[->][blue] (-2.5,4.0)--(7.5,14);
\draw[->][blue] (-2.5,4.5)--(7.0,14);
\draw[->][blue] (-2.5,5.0)--(6.5,14);
\draw[->][blue] (-2.5,5.40)--(6.0,13.90);
\draw[red] (-2.5,5.5)--(-2.5,4.0);
\node at (-2.5,4.0) [fill=red, inner sep=1pt, shape=circle, draw] {};
\node at (-2.5,4.5) [fill=red, inner sep=1pt, shape=circle, draw] {};
\node at (-2.5,5.0) [fill=red, inner sep=1pt, shape=circle, draw] {};
\node at (-2.5,5.5) [fill=red, inner sep=1pt, shape=circle, draw] {};

\draw [->] (-3.0,5.0)-- (6.0,14);
\node at (-3.0,5.0)  [shape = rectangle, draw]{};

\draw[->][blue] (-4.0,4.5)--(5.5,14);
\node at (-4.0,4.5) [fill=red, inner sep=1pt, shape=circle, draw] {};

\draw [->] (-4.0,6.0)-- (4.0,14);
\node at (-4.0,6.0)  [shape = circle, draw]{};

\draw[->][red] (-4.0,5.00)--(5.0,14);
\draw[->][red] (-4.0,5.50)--(4.5,14);
\draw[red] (-4.0,5.0)--(-4.0,5.5);
\node at (-4.0,5.0) [fill=red, inner sep=1pt, shape=circle, draw] {};
\node at (-4.0,5.5) [fill=red, inner sep=1pt, shape=circle, draw] {};

\draw[->][green] (-4.5,4.40)--(5.0,13.90);
\draw[->][green] (-4.5,4.90)--(4.5,13.90);
\draw[->][green] (-4.5,5.40)--(4.0,13.90);
\draw[->][green] (-4.5,6.0)--(3.5,14);
\draw[->][green] (-4.5,6.5)--(3.0,14);
\draw[->][green] (-4.5,7.0)--(2.5,14);
\draw[->][green] (-4.5,7.40)--(2.0,13.90);
\draw[->][green] (-4.5,8.0)--(1.5,14);
\draw[red] (-4.5,7.9)--(-4.5,4.5);
\node at (-4.5,4.5) [fill=red, inner sep=1pt, shape=circle, draw] {};
\node at (-4.5,5.0) [fill=red, inner sep=1pt, shape=circle, draw] {};
\node at (-4.5,5.5) [fill=red, inner sep=1pt, shape=circle, draw] {};
\node at (-4.5,6.0) [fill=red, inner sep=1pt, shape=circle, draw] {};
\node at (-4.5,6.5) [fill=red, inner sep=1pt, shape=circle, draw] {};
\node at (-4.5,7.0) [fill=red, inner sep=1pt, shape=circle, draw] {};
\node at (-4.5,7.5) [fill=red, inner sep=1pt, shape=circle, draw] {};
\node at (-4.5,8.0) [fill=red, inner sep=1pt, shape=circle, draw] {};

\draw [->] (-5.0,7.0)-- (2.0,14);
\node at (-5.0,7.0)  [shape = rectangle, draw]{};

\draw[->][red] (-5.50,7.50)--(1.0,14);
\draw[->][red] (-5.50,8.00)--(0.5,14);
\draw[red] (-5.50,7.50)--(-5.50,8.0);
\node at (-5.50,7.50) [fill=red, inner sep=1pt, shape=circle, draw] {};
\node at (-5.50,8.00) [fill=red, inner sep=1pt, shape=circle, draw] {};

\draw[->][green] (-8.0,5.90)--(0.0,13.90);
\draw[->][green] (-8.0,5.40)--(0.5,13.90);
\draw[->][green] (-8.0,4.90)--(1.0,13.90);
\draw[red] (-8.0,6.0)--(-8.0,5.0);
\node at (-8.0,6.0) [fill=red, inner sep=1pt, shape=circle, draw] {};
\node at (-8.0,5.5) [fill=red, inner sep=1pt, shape=circle, draw] {};
\node at (-8.0,5.0) [fill=red, inner sep=1pt, shape=circle, draw] {};

\draw[->][blue] (-8.0,6.50)--(-0.50,14.00);
\draw[->][blue] (-8.0,7.0)--(-1.00,14);
\draw[->][blue] (-8.0,7.5)--(-1.50,14);
\draw[->][blue] (-8.0,7.90)--(-2.00,13.90);
\draw[red] (-8.0,6.50)--(-8.0,8.0);
\node at (-8.0,6.50) [fill=red, inner sep=1pt, shape=circle, draw] {};
\node at (-8.0,7.0) [fill=red, inner sep=1pt, shape=circle, draw] {};
\node at (-8.0,7.5) [fill=red, inner sep=1pt, shape=circle, draw] {};
\node at (-8.0,8.0) [fill=red, inner sep=1pt, shape=circle, draw] {};

\draw [->] (-6.0,8.0)-- (0.0,14);
\node at (-6.0,8.0)  [shape = circle, draw]{};

\draw [->] (-7.0,9.0)-- (-2.0,14);
\node at (-7.0,9.0)  [shape = rectangle, draw]{};

\draw [->] (-8.0,10.0)-- (-4.0,14);
\node at (-8.0,10.0)  [shape = circle, draw]{};
\draw[->][red] (-8.0,9.5)--(-3.5,14);
\draw[->][red] (-8.0,9.0)--(-3.0,14);
\draw[red] (-8.0,10.0)--(-8.0,9.0);
\draw[->][blue] (-8.0,8.5)--(-2.5,14);
\node at (-8.0,9.5) [fill=red, inner sep=1pt, shape=circle, draw] {};
\node at (-8.0,9.0) [fill=red, inner sep=1pt, shape=circle, draw] {};
\node at (-8.0,8.5) [fill=red, inner sep=1pt, shape=circle, draw] {};

\foreach \y in {1,2,3,4,5,6,7,8,9,10, 11, 12,13,14,15,16,17,18,19}
\node at (-8.5,4.5+\y/2) [fill=red, inner sep=1pt, shape=circle, draw] {};

\draw[red] (-8.5,5.0)--(-8.5,9.50);
\draw[red] (-8.5,10.50)--(-8.5,14.0);

%%%%%%
\draw[->][red] (0,2-0.5)--(12,13.5);
\draw[->][red] (0,2-1.0)--(12,13.0);
\draw[red] (0,2)--(0,-14);

\draw[->][red] (-7,-4.0)--(-12,-9.0);
\draw[->][red] (-7,-4.5)--(-12,-9.5);
\draw[red] (-7,-3.5)--(-7,-4.5);
\node at (-7,-4.0) [fill=red, inner sep=1pt, shape=circle, draw] {};
\node at (-7,-4.5) [fill=red, inner sep=1pt, shape=circle, draw] {};

\draw[->][blue] (0,0.5)--(12,12.5);
\draw[->][blue] (0,-0.10)--(12,11.90);
\draw[->][blue] (0,-0.5)--(12,11.5);
\draw[->][blue] (0,-1.0)--(12,11.0);

\draw[->] (-6,-4.5)--(-12,-10.5);
\node at (-6,-4.5) [shape=circle, draw] {};
\draw[->][blue] (-6,-4.60)--(-12,-10.60);
\draw[->][blue] (-6,-5)--(-12,-11);
\draw[->][blue] (-6,-5.5)--(-12,-11.5);
\draw[->][blue] (-6,-6)--(-12,-12.0);
\draw[red] (-6,-4.55)--(-6,-6.0);
\node at (-6,-4.5) [fill=red, inner sep=1pt, shape=circle, draw] {};
\node at (-6,-5) [fill=red, inner sep=1pt, shape=circle, draw] {};
\node at (-6,-5.5) [fill=red, inner sep=1pt, shape=circle, draw] {};
\node at (-6,-6) [fill=red, inner sep=1pt, shape=circle, draw] {};

\draw[->][green] (0,-1.5)--(12,10.5);
\draw[->][green] (0,-2.10)--(12,9.90);
\draw[->][green] (0,-2.60)--(12,9.40);
\draw[->][green] (0,-3.10)--(12,8.90);
\draw[->][green] (0,-3.5)--(12,8.5);
\draw[->][green] (0,-4.10)--(12,7.90);
\draw[->][green] (0,-4.5)--(12,7.5);
\draw[->][green] (0,-5.0)--(12,7.0);

\draw[->] (-5,-5.5)--(-12,-12.5);
\node at (-5,-5.5) [shape=rectangle, draw] {};
\draw[->][blue] (-4.5,-5.10)--(-12,-12.60);
\node at (-4.5,-5.0) [fill=red, inner sep=1pt, shape=circle, draw] {};
\draw[->][red] (-4.5,-5.50)--(-12,-13.00);
\draw[->][red] (-4.5,-6.00)--(-12,-13.50);
\draw[red] (-4.5,-5.05)--(-4.5,-6.0);
\node at (-4.5,-5.5) [fill=red, inner sep=1pt, shape=circle, draw] {};
\node at (-4.5,-6.0) [fill=red, inner sep=1pt, shape=circle, draw] {};

\draw[->] (-4,-6.5)--(-11.5,-14);
\node at (-4,-6.5) [shape=circle, draw] {};

\draw[->][green] (-4,-5.10)--(-12,-13.10);
\draw[->][green] (-4,-5.60)--(-12,-13.60);
\draw[->][green] (-4,-6)--(-12,-14);
\draw[->][green] (-4,-6.60)--(-11.40,-14.00);
\draw[->][green] (-4,-7)--(-11,-14);
\draw[->][green] (-4,-7.5)--(-10.5,-14);
\draw[->][green] (-4,-8)--(-10,-14);
\draw[->][green] (-4,-8.60)--(-9.40,-14.00);
\draw[red] (-4,-5)--(-4,-8.55);
\node at (-4,-5) [fill=red, inner sep=1pt, shape=circle, draw] {};
\node at (-4,-5.5) [fill=red, inner sep=1pt, shape=circle, draw] {};
\node at (-4,-6.0) [fill=red, inner sep=1pt, shape=circle, draw] {};
\node at (-4,-6.5) [fill=red, inner sep=1pt, shape=circle, draw] {};
\node at (-4,-7.0) [fill=red, inner sep=1pt, shape=circle, draw] {};
\node at (-4,-7.5) [fill=red, inner sep=1pt, shape=circle, draw] {};
\node at (-4,-8.0) [fill=red, inner sep=1pt, shape=circle, draw] {};
\node at (-4,-8.5) [fill=red, inner sep=1pt, shape=circle, draw] {};

\draw[->] (-3,-7.5)--(-9.5,-14);
\node at (-3,-7.5) [shape=rectangle, draw] {};
\draw[->][red] (-3.0,-8.00)--(-9.0,-14.00);
\draw[->][red] (-3.0,-8.50)--(-8.5,-14.00);
\draw[red] (-3.0,-7.50)--(-3.0,-8.50);
\node at (-3.0,-8.0) [fill=red, inner sep=1pt, shape=circle, draw] {};
\node at (-3.0,-8.50) [fill=red, inner sep=1pt, shape=circle, draw] {};

\draw[->][green] (-0.5,-5.60)--(-8.90,-14.00);
\draw[->][green] (-0.5,-6.10)--(-8.40,-14.00);
\draw[->][green] (-0.5,-6.50)--(-8.0,-14);
\draw[red] (-0.5,-5.50)--(-0.5,-6.50);
\node at (-0.5,-5.5) [fill=red, inner sep=1pt, shape=circle, draw] {};
\node at (-0.5,-6.0) [fill=red, inner sep=1pt, shape=circle, draw] {};
\node at (-0.5,-6.5) [fill=red, inner sep=1pt, shape=circle, draw] {};

\draw[->][blue] (-0.5,-7.10)--(-7.40,-14);
\draw[->][blue] (-0.5,-7.50)--(-7.00,-14);
\draw[->][blue] (-0.5,-8.00)--(-6.50,-14);
\draw[->][blue] (-0.5,-8.50)--(-6.0,-14);
\draw[red] (-0.5,-7.00)--(-0.5,-8.50);
\node at (-0.5,-7.0) [fill=red, inner sep=1pt, shape=circle, draw] {};
\node at (-0.5,-7.5) [fill=red, inner sep=1pt, shape=circle, draw] {};
\node at (-0.5,-8.0) [fill=red, inner sep=1pt, shape=circle, draw] {};
\node at (-0.5,-8.5) [fill=red, inner sep=1pt, shape=circle, draw] {};

\draw[->] (-2,-8.5)--(-7.5,-14);
\node at (-2,-8.5) [shape=circle, draw] {};

\draw[->] (-1,-9.5)--(-5.5,-14);
\node at (-1,-9.5) [shape=rectangle, draw] {};
\draw[->][blue] (-0.5,-9.10)--(-5.4,-14.00);
\node at (-0.5,-9.0) [fill=red, inner sep=1pt, shape=circle, draw] {};
\draw[->][red] (-0.5,-9.5)--(-5.0,-14);
\draw[->][red] (-0.5,-10.0)--(-4.5,-14);
\draw[red] (-0.5,-9)--(-0.5,-10.0);
\node at (-0.5,-9.5) [fill=red, inner sep=1pt, shape=circle, draw] {};
\node at (-0.5,-10.0) [fill=red, inner sep=1pt, shape=circle, draw] {};

\draw[->] (0,-10.5)--(-3.5,-14);
\node at (0,-10.5) [shape=circle, draw] {};

%\foreach \y in {15,16,17}
%\draw (0.1,2.1-\y/2) node[anchor=east,red] {\tiny{$d_2$}};
\draw (0.15,-5.3) node[anchor=east,red] {\tiny{$d_2$}};
\draw[->,red] (0,-5.5) to [ bend left=55] (-0.5,-5.60) ;
\draw (0.15,-5.8) node[anchor=east,red] {\tiny{$d_2$}};
\draw[->,red] (0,-6.0) to [ bend left=55] (-0.5,-6.10) ;
\draw (0.15,-6.3) node[anchor=east,red] {\tiny{$d_2$}};
\draw[->,red] (0,-6.5) to [ bend left=55] (-0.5,-6.60) ;
\draw (0.15,-8.8) node[anchor=east,red] {\tiny{$d_2$}};
\draw[->,red] (0,-9.0) to [ bend left=55] (-0.5,-9.10) ;
\draw (0.15,-9.3) node[anchor=east,red] {\tiny{$d_3$}};
\draw[->,red] (0,-9.5) to [ bend left=55] (-0.5,-9.60) ;
\draw (0.15,-9.8) node[anchor=east,red] {\tiny{$d_3$}};
\draw[->,red] (0,-10.0) to [ bend left=55] (-0.5,-10.10) ;
\draw (0.15,-6.8) node[anchor=east,red] {\tiny{$d_2$}};
\draw[->,red] (0,-7.00) to [ bend left=55] (-0.5,-7.10) ;
\draw (0.15,-7.3) node[anchor=east,red] {\tiny{$d_2$}};
\draw[->,red] (0,-7.50) to [ bend left=55] (-0.5,-7.60) ;
\draw (0.15,-7.8) node[anchor=east,red] {\tiny{$d_2$}};
\draw[->,red] (0,-8.0) to [ bend left=55] (-0.5,-8.10) ;
\draw (0.15,-8.3) node[anchor=east,red] {\tiny{$d_2$}};
\draw[->,red] (0,-8.5) to [ bend left=55] (-0.5,-8.60) ;

\draw (-3.83,-4.9) node[anchor=east,red] {\tiny{$d_2$}};
\draw[->,red] (-4,-5.0) to [ bend left=55] (-4.5,-5.10) ;
\draw (-3.83,-5.4) node[anchor=east,red] {\tiny{$d_2$}};
\draw[->,red] (-4,-5.5) to [ bend left=55] (-4.5,-5.60) ;
\draw (-3.83,-5.9) node[anchor=east,red] {\tiny{$d_2$}};
\draw[->,red] (-4,-6.0) to [ bend left=55] (-4.5,-6.10) ;

\draw [->](1,1)-- (12,12);

\node at (1,1) [shape=circle, draw] {};
\draw (1,1) node[anchor=east]{$2u$};

\foreach \y in {1,2,3,4,5,6,7,8,9,10,11,12,13,14,15,16,17,18,19,20,21,22}
\node at (1+\y/2,1+\y/2) [shape=circle, draw] {};

\draw [->](2,0)-- (12,10);
\node at (2,0) [shape=circle, draw] {};

\draw (2,0) node[anchor=east]{\small{$2u^2$}};
\node at (2.5,0.5) [shape=rectangle, draw] {};
\draw (2.5,0.5) node[anchor=west]{$u^2 \bar v_1$};

\draw[->][red] (2.5,0)--(12,9.5);
\draw[->][red] (2.5,-0.5)--(12,9.0);
\draw[red] (2.5,0.5)--(2.5,-0.5);
\node at (2.5,0) [fill=red, inner sep=1pt, shape=circle, draw] {};
\node at (2.5,-0.5) [fill=red, inner sep=1pt, shape=circle, draw] {};

\draw [->](3,-1)-- (12,8.0);
\node at (3,-1) [shape=circle, draw] {};

\draw (3,-1) node[anchor=east]{\small{$2u^3$}};

\draw [->](4,-2)-- (12,6.0);
\node at (4,-2) [shape=circle, draw] {};

\draw (4,-2) node[anchor=east]{\small{$2u^4$}};
\node at (4.5,-1.5) [shape=rectangle, draw] {};
\node at (5.5,-0.5) [fill=blue, inner sep=3pt, shape=rectangle, draw] {};

\draw (4.5,-1.5) node[anchor=east]{\small{$u^4 \bar v_1$}};
\draw (5.5,-0.5) node[anchor=east]{\small{$u^4 \bar v_2$}};

\draw[->][red] (4.5,-2.0)--(12,5.5);
\draw[->][red] (4.5,-2.5)--(12,5.0);
%\draw[->][blue] (5.5,-1.0)--(12,5.5);
%\draw[->][blue] (5.5,-1.5)--(12,5.0);
\draw[->][blue] (5.5,-2.0)--(12,4.5);
\draw[->][blue] (5.5,-2.60)--(12,3.90);
\draw[->][blue] (5.5,-3.0)--(12,3.5);
\draw[->][blue] (5.5,-3.5)--(12,3.0);
\draw[red] (4.5,-1.5)--(4.5,-2.5);
\draw[red] (5.5,-0.5)--(5.5,-3.5);
\node at (4.5,-2.0) [fill=red, inner sep=1pt, shape=circle, draw] {};
\node at (4.5,-2.5) [fill=red, inner sep=1pt, shape=circle, draw] {};
\node at (5.5,-1.0) [fill=red, inner sep=1pt, shape=circle, draw] {};
\node at (5.5,-1.5) [fill=red, inner sep=1pt, shape=circle, draw] {};
\node at (5.5,-2.0) [fill=red, inner sep=1pt, shape=circle, draw] {};
\node at (5.5,-2.5) [fill=red, inner sep=1pt, shape=circle, draw] {};
\node at (5.5,-3.0) [fill=red, inner sep=1pt, shape=circle, draw] {};
\node at (5.5,-3.5) [fill=red, inner sep=1pt, shape=circle, draw] {};

\draw [->](5,-3)-- (12,4.0);

\node at (5,-3) [shape=circle, draw] {};
\draw (5,-3) node[anchor=east]{\small{$2u^5$}};

\draw [->](6,-4)-- (12,2.0);
\node at (6,-4) [shape=circle, draw] {};

\draw (6,-4) node[anchor=east]{\small{$2u^6$}};
\node at (6.5,-3.5) [shape=rectangle, draw] {};
\draw (6.5,-3.5) node[anchor=west]{\small{$u^6 \bar v_1$}};

\draw[->][red] (6.5,-4.0)--(12,1.5);
\draw[->][red] (6.5,-4.5)--(12,1.0);
\draw[red] (6.5,-3.5)--(6.5,-4.5);
\node at (6.5,-4.0) [fill=red, inner sep=1pt, shape=circle, draw] {};
\node at (6.5,-4.5) [fill=red, inner sep=1pt, shape=circle, draw] {};

\draw [->](7,-5)-- (12,0.0);
\node at (7,-5) [shape=circle, draw] {};
\draw (7,-5) node[anchor=east]{\small{$2u^7$}};

\node at (-4.0,-0.25) [fill=yellow, inner sep=1.5pt, shape=diamond, draw] {};

\end{tikzpicture}$$
\caption{The basic block, SS(GBB) and $\mathbb{Z}^{\rm{SS(GBB)}}$ on $BP\mathbb R\langle 3 \rangle^{Q}_{\rost}$} \label{fig:GL1}
\end{figure}

Notice that connectivity is the condition that at $x+y\sigma$, when
$x<0$ then there are no entries with $x+y<0$. This forces all the differentials $d_2$ and $d_3$
as indicated. It also forces the additive extension at $0-25\sigma$
and $-8-17\sigma$.

\newpage

Turning now to the negative block $NB$, Tables \ref{tab:table4} and \ref{tab:table5} are the same as
Table \ref{tab:table1}, except that in the $u^0$ column, each module $M$ is replaced
by $\ker (M\lra \Ftwo)$ for $\delta \geq 0$. There is also a column of
$\Ftwo$s starting at $(-1+\sigma)u^0$ and going upwards. Since this is
evidently all $J$-torsion,  we have omitted it from Tables 4 and 6. 
%\ref{tab:table4, tab:table6}.
\renewcommand{\tabcolsep}{14.00pt}
\renewcommand{\arraystretch}{0.70}
\newcolumntype{a}{>{\columncolor{Gray}}c}
%\newcolumntype{b}{>{\columncolor{white}}c}
\begin{table}[h!]
\begin{center}
 \caption{Negative block of $BP \mathbb R\langle 3
   \rangle^{Q}_{\rost}$ for $1, u, u^2 , u^3$ and $\delta = 0, 1,2,3, \ldots, 15$}
    \label{tab:table4}
    \begin{tabular}{ ||!{\mystrut}a | p{2.50cm} | p{2.90cm} | p{2.90cm} |p{2.50cm}  || }
    \hline
      $\delta$ &\cellcolor{yellow!55}$1$ &\cellcolor{yellow!55}$u$ &\cellcolor{yellow!55}$u^2$ &\cellcolor{yellow!55}$u^3$  \\ [1.5ex]
    \hline \hline
      $0$ &$(2,\bar v_1,\bar v_2, \bar v_3)P$        &       &            &       \\    \hline
      $1$ &$(\bar v_1,\bar v_2,\bar v_3)\Pb_0$  &       &            &       \\    \hline
      $2$ &$(\bar v_1,\bar v_2, \bar v_3)\Pb_0$  &       &            &       \\    \hline
      $3$ &$(\bar v_2,\bar v_3)\Pb_1$      &       &            &       \\    \hline
      $4$ &$(\bar v_2,\bar v_3)\Pb_1$      &$(2)P$ &            &      \\    \hline
      $5$ &$(\bar v_2,\bar v_3)\Pb_1$      &       &            &       \\    \hline
      $6$ &$(\bar v_2,\bar v_3)\Pb_1$      &       &            &       \\    \hline
      $7$ &$(\bar v_3)\Pb_2$      &   &       &        \\    \hline
      $8$ &$(\bar v_3)\Pb_2$  &       &$(2,\bar v_1)P$  &       \\    \hline
      $9$ &$(\bar v_3)\Pb_2$  &       &$(\bar v_1)\Pb_0$  &       \\    \hline
      $10$ &$(\bar v_3)\Pb_2$ &       &$(\bar v_1)\Pb_0$  &       \\    \hline
      $11$ &$(\bar v_3)\Pb_2$ &       &                      &       \\    \hline
      $12$ &$(\bar v_3)\Pb_2$ &       & &$(2)P$           \\ \hline
      $13$ &$(\bar v_3)\Pb_2$ &       & &           \\ \hline
      $14$ &$(\bar v_3)\Pb_2$ &       & &           \\ \hline
      $15$ & &       &            &       \\
    \hline
    \end{tabular}
\end{center}
\end{table}
\vspace{0.1cm}

\renewcommand{\tabcolsep}{14.00pt}
\renewcommand{\arraystretch}{0.70}
\newcolumntype{a}{>{\columncolor{Gray}}c}
%\newcolumntype{b}{>{\columncolor{white}}c}
\begin{table}[h!]
\begin{center}
 \caption{Negative block of $BP \mathbb R\langle 3 \rangle^{Q}_{\rost}$ for $u^4, u^5, u^6, u^7$ and $\delta = 16,17,18, \ldots, 31$}
    \label{tab:table5}
    \begin{tabular}{ ||!{\mystrut}a | p{2.50cm}|p{2.90cm} |p{2.90cm} |p{2.50cm} || }
    \hline
      $\delta$ &\cellcolor{yellow!55}$u^4$
      &\cellcolor{yellow!55}$u^5$ &\cellcolor{yellow!55}$u^6$ &\cellcolor{yellow!55}$u^7$  \\ [1.5ex]
    \hline \hline
      $16$ &$(2,\bar v_1,\bar v_2)P$            & & &\\ \hline
      $17$ &$(\bar v_1,\bar v_2)\Pb_0$       & & &\\ \hline
      $18$ &$(\bar v_1,\bar v_2)\Pb_0$       & & &\\ \hline
      $19$ &$(\bar v_2)\Pb_1$           & & &\\ \hline
      $20$        &$(\bar v_2)\Pb_1$&$(2)P$   & &\\ \hline
      $21$        &$(\bar v_2)\Pb_1$&   & &\\ \hline
      $22$        &$(\bar v_2)\Pb_1$&   & &\\ \hline
      $23$        &                    &       & &\\ \hline
      $24$        &                    &&$(2,\bar v_1)P$      & \\ \hline
      $25$        &                    &&$(\bar v_1)\Pb_0$ & \\ \hline
      $26$        &                    &&$(\bar v_1)\Pb_0$ & \\ \hline
      $27$        &                    & &           &\\ \hline
      $28$        &                    & & &$(2)P$    \\ \hline
      $29$        &                    & &           &\\ \hline
      $30$        &                    & &           &\\ \hline
      $31$        &                    & &           &\\
    \hline
    \end{tabular}
\end{center}
\end{table}
\vspace{0.1cm}

%\newpage

\renewcommand{\tabcolsep}{14.00pt}
\renewcommand{\arraystretch}{0.80}
\newcolumntype{a}{>{\columncolor{Gray}}c}
%\newcolumntype{b}{>{\columncolor{white}}c}
\begin{table}[h!]
\begin{center}
 \caption{Local cohomology of the negative block for $1, u, u^2$ and $u^3$}
    \label{tab:table6}
    \begin{tabular}{ ||!{\mystrut}a | p{3.10cm} | p{1.80cm} | p{4.10cm} |p{1.80cm}  || }
    \hline
      $\delta$ &\cellcolor{yellow!55}$1$ &\cellcolor{yellow!55}$u$ &\cellcolor{yellow!55}$u^2$ &\cellcolor{yellow!55}$u^3$  \\ [1.5ex]
    \hline \hline
      $-3$ &\textcolor{red}{$P^*(-14\rho)$}        &       &            &       \\    \hline
      $-2$ &\textcolor{red}{$\Pb_0^\vee(-14\rho)$}  &       &            &       \\    \hline
      $-1$ &\textcolor{red}{$\Pb_0^\vee(-14\rho)$} $\oplus$ $\mathbb{F}_2$ &       &            &       \\    \hline
      $0$ &\hspace*{12.5ex} $\mathbb{F}_2$                                                &       &            &       \\    \hline
      $1$ &\textcolor{blue}{$\Pb_1^\vee(-12\rho)$} $\oplus$ $\mathbb{F}_2$    &\textcolor{red}{$P^*(-14\rho)$} &            &      \\    \hline
      $2$ &\textcolor{blue}{$\Pb_1^\vee(-12\rho)$} $\oplus$ $\mathbb{F}_2$      &       &            &       \\    \hline
      $3$ &\textcolor{blue}{$\Pb_1^\vee(-12\rho)$} $\oplus$ $\mathbb{F}_2$      &       &            &       \\    \hline
      $4$ &\textcolor{blue}{$\Pb_1^\vee(-12\rho)$} $\oplus$ $\mathbb{F}_2$      &   &       &        \\    \hline
      $5$ &\hspace*{13.5ex}$\mathbb{F}_2$  &       &\textcolor{red}{$P^*(-14\rho)$} $\oplus$ \textcolor{blue}{$\Pb_1^\vee(-13\rho)$}  &       \\    \hline
      $6$ &\textcolor{green}{$\Pb_2^\vee(-8\rho), d_2$}  \hspace*{-0.5ex}$\oplus$ $\mathbb{F}_2$ &       &\textcolor{red}{$\Pb_0^\vee(-13\rho)$}  &       \\    \hline
      $7$ &\textcolor{green}{$\Pb_2^\vee(-8\rho), d_2$} \hspace*{-0.5ex}$\oplus$ $\mathbb{F}_2$  &       &\textcolor{red}{$\Pb_0^\vee(-13\rho)$}  &       \\    \hline
      $8$ &\textcolor{green}{$\Pb_2^\vee(-8\rho), d_2$} \hspace*{-0.5ex}$\oplus$ $\mathbb{F}_2$ &       &                      &       \\    \hline
      $9$ &\textcolor{green}{$\Pb_2^\vee(-8\rho)$} $\oplus$ $\mathbb{F}_2$  &       & &\textcolor{red}{$P^*(-14\rho)$}          \\ \hline
      $10$ &\textcolor{green}{$\Pb_2^\vee(-8\rho)$} $\oplus$ $\mathbb{F}_2$  &       & &           \\ \hline
      $11$ &\textcolor{green}{$\Pb_2^\vee(-8\rho)$} $\oplus$ $\mathbb{F}_2$  &       & &           \\ \hline
      $12$ &\textcolor{green}{$\Pb_2^\vee(-8\rho)$} $\oplus$ $\mathbb{F}_2$ &       &            &       \\ \hline
      $13$ &\textcolor{green}{$\Pb_2^\vee(-8\rho)$} $\oplus$ $\mathbb{F}_2$ &       &            &       \\
    \hline
    \end{tabular}
\end{center}
\end{table}
\vspace{0.2cm}

\renewcommand{\tabcolsep}{14.00pt}
\renewcommand{\arraystretch}{0.90}
\newcolumntype{a}{>{\columncolor{Gray}}c}
%\newcolumntype{b}{>{\columncolor{white}}c}
\begin{table}[h!]
\begin{center}
 \caption{Local cohomology of the negative block for $u^4, u^5, u^6$ and $u^7$}
    \label{tab:table7}
    \begin{tabular}{ ||!{\mystrut}a | p{4.00cm}|p{1.40cm} |p{4.20cm} |p{1.40cm} || }
    \hline
      $\delta$ &\cellcolor{yellow!55}$u^4$
      &\cellcolor{yellow!55}$u^5$ &\cellcolor{yellow!55}$u^6$ &\cellcolor{yellow!55}$u^7$  \\ [1.5ex]
    \hline \hline
      $13$ &\textcolor{red}{$P^*(-14\rho)$}    & & &\\ \hline
      $14$ &\textcolor{red}{$\Pb_0^\vee(-14\rho)$} $\oplus$ \textcolor{blue}{$\Pb_2^\vee(-9\rho)$}     & & &\\ \hline
      $15$ &\textcolor{red}{$\Pb_0^\vee(-14\rho)$} $\oplus$ \textcolor{blue}{$\Pb_2^\vee(-9\rho)$}      & & &\\ \hline
      $16$ &\hspace*{13ex}\textcolor{blue}{$\Pb_2^\vee(-9\rho)$}           & & &\\ \hline
      $17$        &\textcolor{blue}{$\Pb_1^\vee(-10\rho)$}   &\textcolor{red}{$P^*(-14\rho)$}   & &\\ \hline
      $18$        &\textcolor{blue}{$\Pb_1^\vee(-10\rho)$}   &   & &\\ \hline
      $19$        &\textcolor{blue}{$\Pb_1^\vee(-10\rho)$}   &   & &\\ \hline
      $20$        &\textcolor{blue}{$\Pb_1^\vee(-10\rho)$}      &   & &\\ \hline
      $21$        &                    &&\textcolor{red}{$P^*(-14\rho)$} $\oplus$ \textcolor{blue}{$\Pb_1^\vee(-13\rho)$}      & \\ \hline
      $22$        &                    &&\textcolor{red}{$\Pb_0^\vee(-13\rho)$}  & \\ \hline
      $23$        &                    &&\textcolor{red}{$\Pb_0^\vee(-13\rho)$} & \\ \hline
      $24$        &                    & &           &\\ \hline
      $25$        &                    & & &\textcolor{red}{$P^*(-14\rho)$}    \\ \hline
      $26$        &                    & &           &\\
    \hline
    \end{tabular}
\end{center}
\end{table}
\vspace{0.2cm}
As in the case of the local cohomology of the basic block, we have
coloured $H^1_J$-groups in green, $H^2_J$-groups in blue and $H^3_J$-groups in
red, and with the same caveat.

%\newpage

\begin{figure}
\centering
$$\begin{tikzpicture}[scale =0.688]
%\clip (-9, -8.6) rectangle (5, 8.5);
\draw[step=0.5, gray, very thin] (-12,-14) grid (12, 14);
\draw (-0.5,5.0 ) node[anchor=east, draw=orange]{\Large{NB}};

\foreach \y in {1,2,3,4,5,6,7,8}
\node at (-4.5,10+\y/2) [fill=red, inner sep=1pt, shape=circle, draw] {};
\draw[red] (-4.5,10.50)--(-4.5,14.00);

%%%%%%(NB)
\draw [->] (3,3)-- (12,12);
\node at (3,3)  [shape = circle, draw]{};

\draw[->][red] (2.5,3.5)--(12.0,13.0);
\draw[->][red] (2.5,4.0)--(12.0,13.5);
\draw[red] (2.5,3.5)--(2.5,4.5);
\node at (2.5,3.5) [fill=red, inner sep=1pt, shape=circle, draw] {};
\node at (2.5,4.0) [fill=red, inner sep=1pt, shape=circle, draw] {};

\draw [->] (2.0,4.0)-- (12.0,14);
\node at (2.0,4.0)  [shape = circle, draw]{};
\node at (2.5,4.5) [shape=rectangle, draw] {};
\node at (0.5,6.5) [shape=rectangle, draw] {};
\node at (1.5,7.5) [fill=blue, inner sep=3pt, shape=rectangle, draw] {};
\node at (-3.5,10.5) [shape=rectangle, draw] {};
\node at (-2.5,11.5) [fill=blue, inner sep=3pt, shape=rectangle, draw] {};
\node at (-0.5,13.5) [fill=green, inner sep=3pt, shape=rectangle, draw] {};
\node at (-1.5,8.5) [shape=rectangle, draw] {};

\draw (-3.5,10.5) node[anchor=east]{\small{$\bar v_1$}};
\draw (-2.5,11.5) node[anchor=east]{\small{$\bar v_2$}};
\draw (-0.5,13.5) node[anchor=east]{\small{$\bar v_3$}};

\draw (0.5,6.5) node[anchor=east]{\small{$\bar v_1$}};
\draw (1.5,7.5) node[anchor=east]{\small{$\bar v_2$}};
\draw (2.5,4.5) node[anchor=east]{\small{$\bar v_1$}};
\draw (-1.5,8.5) node[anchor=east]{\small{$\bar v_1$}};

\draw[->][blue] (1.5,4.5)--(11.0,14);
\draw[->][blue] (1.5,5.0)--(10.5,14);
\draw[->][blue] (1.5,5.40)--(10.0,13.90);
\draw[->][blue] (1.5,6.0)--(9.5,14);
\draw[red] (1.5,7.5)--(1.5,4.5);
\node at (1.5,4.5) [fill=red, inner sep=1pt, shape=circle, draw] {};
\node at (1.5,5.0) [fill=red, inner sep=1pt, shape=circle, draw] {};
\node at (1.5,5.5) [fill=red, inner sep=1pt, shape=circle, draw] {};
\node at (1.5,6.0) [fill=red, inner sep=1pt, shape=circle, draw] {};
\node at (1.5,6.5) [fill=red, inner sep=1pt, shape=circle, draw] {};
\node at (1.5,7.0) [fill=red, inner sep=1pt, shape=circle, draw] {};

\draw [->] (1.0,5.0)-- (10.0,14);
\node at (1.0,5.0)  [shape = circle, draw]{};

\draw [->] (0.0,6.0)-- (8.0,14);
\node at (0.0,6.0)  [shape = circle, draw]{};
\draw[->][red] (0.5,6.0)--(8.5,14);
\draw[->][red] (0.5,5.5)--(9.0,14);
\draw[red] (0.5,6.5)--(0.5,5.5);
\node at (0.5,6.0) [fill=red, inner sep=1pt, shape=circle, draw] {};
\node at (0.5,5.5) [fill=red, inner sep=1pt, shape=circle, draw] {};

\draw[->][green] (-0.5,6.5)--(7.0,14);
\draw[->][green] (-0.5,7.0)--(6.5,14);
\draw[->][green] (-0.5,7.40)--(6.0,13.90);
\draw[->][green] (-0.5,8.00)--(5.5,14.00);
\draw[->][green] (-0.5,8.40)--(5.0,13.90);
\draw[->][green] (-0.5,8.90)--(4.5,13.90);
\draw[->][green] (-0.5,9.40)--(4.0,13.90);
\draw[->][green] (-0.5,10.0)--(3.5,14);
\draw[red] (-0.5,10.0)--(-0.5,6.5);
\node at (-0.5,6.5) [fill=red, inner sep=1pt, shape=circle, draw] {};
\node at (-0.5,7.0) [fill=red, inner sep=1pt, shape=circle, draw] {};
\node at (-0.5,7.5) [fill=red, inner sep=1pt, shape=circle, draw] {};
\node at (-0.5,8.0) [fill=red, inner sep=1pt, shape=circle, draw] {};
\node at (-0.5,8.5) [fill=red, inner sep=1pt, shape=circle, draw] {};
\node at (-0.5,9.0) [fill=red, inner sep=1pt, shape=circle, draw] {};
\node at (-0.5,9.5) [fill=red, inner sep=1pt, shape=circle, draw] {};
\node at (-0.5,10.0) [fill=red, inner sep=1pt, shape=circle, draw] {};
\node at (-0.5,10.5) [fill=red, inner sep=1pt, shape=circle, draw] {};
\node at (-0.5,11.0) [fill=red, inner sep=1pt, shape=circle, draw] {};
\node at (-0.5,11.5) [fill=red, inner sep=1pt, shape=circle, draw] {};
\node at (-0.5,12.0) [fill=red, inner sep=1pt, shape=circle, draw] {};
\node at (-0.5,12.5) [fill=red, inner sep=1pt, shape=circle, draw] {};
\node at (-0.5,13.0) [fill=red, inner sep=1pt, shape=circle, draw] {};

\draw [->] (-1.0,7.0)-- (6.0,14);
\node at (-1.0,7.0)  [shape = circle, draw]{};

\draw [->] (-2.0,8.0)-- (4.0,14);
\node at (-2.0,8.0)  [shape = circle, draw]{};

\draw [->] (-3.0,9.0)-- (2.0,14);
\node at (-3.0,9.0)  [shape = circle, draw]{};

\draw[->][red] (-1.5,8.0)--(4.5,14);
\draw[->][red] (-1.5,7.5)--(5.0,14);
\draw[red] (-1.5,7.5)--(-1.5,8.5);
\node at (-1.5,8.0) [fill=red, inner sep=1pt, shape=circle, draw] {};
\node at (-1.5,7.5) [fill=red, inner sep=1pt, shape=circle, draw] {};

\draw[->][blue] (-2.5,8.5)--(3.0,14);
\draw[->][blue] (-2.5,9.0)--(2.5,14);
\draw[->][blue] (-2.5,9.40)--(2.0,13.90);
\draw[->][blue] (-2.5,10.0)--(1.5,14);
\draw[red] (-2.5,11.5)--(-2.5,8.5);
\node at (-2.5,8.5) [fill=red, inner sep=1pt, shape=circle, draw] {};
\node at (-2.5,9.0) [fill=red, inner sep=1pt, shape=circle, draw] {};
\node at (-2.5,9.5) [fill=red, inner sep=1pt, shape=circle, draw] {};
\node at (-2.5,10.0) [fill=red, inner sep=1pt, shape=circle, draw] {};
\node at (-2.5,10.5) [fill=red, inner sep=1pt, shape=circle, draw] {};
\node at (-2.5,11.0) [fill=red, inner sep=1pt, shape=circle, draw] {};

\draw [->] (-4.0,10.0)-- (0.0,14);
\node at (-4.0,10.0)  [shape = circle, draw]{};
\draw[->][red] (-3.5,10.0)--(0.5,14);
\draw[->][red] (-3.5,9.5)--(1.0,14);
\draw[red] (-3.5,10.5)--(-3.5,9.5);
\node at (-3.5,10.0) [fill=red, inner sep=1pt, shape=circle, draw] {};
\node at (-3.5,9.5) [fill=red, inner sep=1pt, shape=circle, draw] {};

\draw[->][red, very thin] (-4.0,10.0)--(3,3);

%%%%%SS(SS(GNB))
\draw (-6.5,6.5 ) node[anchor=east, draw=orange]{\Large{SS(GNB)}};

\foreach \y in {1,2,3,4,5,6,7,8,9,10,11,12,13,14,15,16,17,18,19,20,21,22,23,24,25}
\node at (-4.5,10.5-\y/2) [fill=red, inner sep=1pt, shape=circle, draw] {};
\draw[red] (-4.5,10.0)--(-4.5,-0.50);
\draw[red] (-4.5,-1.00)--(-4.5,-2.00);
\draw[line width=0.20mm, orange] (-4.5,10.0)--(-4.5,10.50);
\draw[line width=0.20mm, orange] (-4.5,-0.50)--(-4.5,-1.00);

\draw[orange] (-4.5,10.0) to [bend left=55] (-4.5,10.50) ;
\draw[orange] (-4.5,-0.50) to [bend left=55] (-4.5,-1.00) ;

\draw [->] (-11,4.5)-- (-12,3.5);
\node at (-11,4.5)  [shape = circle, draw]{};

\draw[->][red] (-11,4.0)--(-12,3.0);
\draw[->][red] (-11,3.5)--(-12,2.5);
\draw[red] (-11,4.5)--(-11,3.5);
\node at (-11,4.0) [fill=red, inner sep=1pt, shape=circle, draw] {};
\node at (-11,3.5) [fill=red, inner sep=1pt, shape=circle, draw] {};

\draw [->] (-10.0,3.5)-- (-12.0,1.5);
\node at (-10.0,3.5)  [shape = circle, draw]{};

\draw[->][blue] (-8.5-1.5,4.90-1.5)--(-12.0,1.40);
\draw[->][blue] (-8.5-1.5,4.5-1.5)--(-12.0,1.0);
\draw[->][blue] (-8.5-1.5,4.0-1.5)--(-12.0,0.5);
\draw[->][blue] (-8.5-1.5,3.5-1.5)--(-12.0,0.0);
\draw[red] (-8.5-1.5,5.0-1.5)--(-8.5-1.5,3.5-1.5);
\node at (-8.5-1.5,5.0-1.6) [fill=red, inner sep=1pt, shape=circle, draw] {};
\node at (-8.5-1.5,4.5-1.5) [fill=red, inner sep=1pt, shape=circle, draw] {};
\node at (-8.5-1.5,4.0-1.5) [fill=red, inner sep=1pt, shape=circle, draw] {};
\node at (-8.5-1.5,3.5-1.5) [fill=red, inner sep=1pt, shape=circle, draw] {};

\draw [->] (-9.0,2.5)-- (-12.0,-0.5);
\node at (-9.0,2.5)  [shape = circle, draw]{};
\draw[->][blue] (-8.5,2.90)--(-12.0,-0.60);
\node at (-8.5,3.0) [fill=red, inner sep=1pt, shape=circle, draw] {};

\draw[->][red] (-8.5,2.5)--(-12,-1.0);
\draw[->][red] (-8.5,2.0)--(-12,-1.5);
\draw[red] (-8.5,2.5)--(-8.5,2.0);
\node at (-8.5,2.5) [fill=red, inner sep=1pt, shape=circle, draw] {};
\node at (-8.5,2.0) [fill=red, inner sep=1pt, shape=circle, draw] {};

\draw[->][green] (-4.5-3.5,6.40-3.5)--(-12,-1.1);
\draw[->][green] (-4.5-3.5,5.90-3.5)--(-12,-1.6);
\draw[->][green] (-4.5-3.5,5.50-3.5)--(-12,-2.00);
\draw[->][green] (-4.5-3.5,4.90-3.5)--(-12,-2.60);
\draw[->][green] (-4.5-3.5,4.50-3.5)--(-12,-3.00);
\draw[->][green] (-4.5-3.5,4.00-3.5)--(-12,-3.50);
\draw[->][green] (-4.5-3.5,3.50-3.5)--(-12,-4.00);
\draw[->][green] (-4.5-3.5,2.90-3.5)--(-12,-4.60);
\draw[->][green] (-4.5,2.40)--(-12,-5.10);
\draw[->][green] (-4.5,1.90)--(-12,-5.60);
\draw[->][green] (-4.5,1.50)--(-12,-6.00);
\draw[red] (-4.5-3.5,6.40-3.5)--(-4.5-3.5,2.90-3.5);
\node at (-4.5-3.5,6.40-3.5) [fill=red, inner sep=1pt, shape=circle, draw] {};
\node at (-4.5-3.5,5.90-3.5) [fill=red, inner sep=1pt, shape=circle, draw] {};
\node at (-4.5-3.5,5.50-3.5) [fill=red, inner sep=1pt, shape=circle, draw] {};
\node at (-4.5-3.5,4.90-3.5) [fill=red, inner sep=1pt, shape=circle, draw] {};
\node at (-4.5-3.5,4.50-3.5) [fill=red, inner sep=1pt, shape=circle, draw] {};
\node at (-4.5-3.5,4.00-3.5) [fill=red, inner sep=1pt, shape=circle, draw] {};
\node at (-4.5-3.5,3.50-3.5) [fill=red, inner sep=1pt, shape=circle, draw] {};
\node at (-4.5-3.5,2.90-3.5) [fill=red, inner sep=1pt, shape=circle, draw] {};

\draw (-4.5-3.5+0.2,6.40-3.5+0.2) node[anchor=east,red] {\tiny{$d_2$}};
\draw[->,red] (-4.5-3.5,6.40-3.5) to [ bend left=55] (-4.5-3.5-0.5,6.40-3.5) ;
\draw (-4.5-3.5+0.2,5.90-3.5+0.2) node[anchor=east,red] {\tiny{$d_2$}};
\draw[->,red] (-4.5-3.5,5.90-3.5) to [ bend left=55] (-4.5-3.5-0.5,5.90-3.5) ;
\draw (-4.5-3.5+0.2,5.40-3.5+0.2) node[anchor=east,red] {\tiny{$d_2$}};
\draw[->,red] (-4.5-3.5,5.50-3.5) to [ bend left=55] (-4.5-3.5-0.5,5.40-3.5) ;

\draw [->] (-8.0,1.5)-- (-12,-2.5);
\node at (-8.0,1.5)  [shape = circle, draw]{};
\draw [->] (-7.0,0.5)-- (-12,-4.5);
\node at (-7.0,0.5)  [shape = circle, draw]{};
\draw [->] (-6.0,-0.50)-- (-12.0,-6.5);
\node at (-6.0,-0.50)  [shape = circle, draw]{};
\draw [->] (-5.0,-1.5)-- (-12,-8.5);
\node at (-5.0,-1.5)  [shape = circle, draw]{};
\draw [->] (-4.0,-2.50)-- (-12.0,-10.5);
\node at (-4.0,-2.50)  [shape = circle, draw]{};

\draw[->][red] (-7.0,0.0)--(-12,-5.0);
\draw[->][red] (-7.0,-0.5)--(-12,-5.5);
\draw[red] (-7.0,0.50)--(-7.0,-0.5);
\node at (-7.0,0.0) [fill=red, inner sep=1pt, shape=circle, draw] {};
\node at (-7.0,-0.5) [fill=red, inner sep=1pt, shape=circle, draw] {};

\draw[->][blue] (-4.5,0.9)--(-12,-6.6);
\draw[->][blue] (-4.5,0.5)--(-12,-7.0);
\draw[->][blue] (-4.5,0.0)--(-12,-7.5);
\draw[->][blue] (-4.5,-0.5)--(-12,-8.0);

\draw[->][blue] (-4.5,-1.1)--(-12,-8.6);
\draw[->][red] (-4.5,-1.5)--(-12,-9.0);
\draw[->][red] (-4.5,-2.0)--(-12,-9.5);

%%%%%SS({$\mathbb Z}^SSGNB)
\draw (3.5,-5.0 ) node[anchor=east, draw=orange]{\Large{${\mathbb Z}^{\rm{SS(GNB)}}$}};

\foreach \y in {1,2,3,4,5,6,7,8,9,10,11,12,13,14,15,16,17,18,19,20,21,22,23,24,25,26,27,28,29,30,31,32}
\node at (4.0,2-\y/2) [fill=red, inner sep=1pt, shape=circle, draw] {};
\draw[red] (4.0,2.0)--(4.0,-14.00);

\draw [->] (11,-5)-- (12,-4);
\node at (11,-5)  [shape = circle, draw]{};

\draw[->][red] (10.5,-4.5)--(12,-3.0);
\draw[->][red] (10.5,-4.0)--(12,-2.5);
\draw[red] (10.5,-4.5)--(10.5,-4.0);
\node at (10.5,-4.5) [fill=red, inner sep=1pt, shape=circle, draw] {};
\node at (10.5,-4.0) [fill=red, inner sep=1pt, shape=circle, draw] {};

\draw [->] (10.0,-4.0)-- (12.0,-2.0);
\node at (10.0,-4.0)  [shape = circle, draw]{};

\draw[->][blue] (8.0+1.5,-4.10+1.5)--(12.0,-0.10);
\draw[->][blue] (8.0+1.5,-4.5+1.5)--(12.0,-0.5);
\draw[->][blue] (8.0+1.5,-5.0+1.5)--(12.0,-1.0);
\draw[->][blue] (8.0+1.5,-5.5+1.5)--(12.0,-1.5);
\draw[red] (8.0+1.5,-4.0+1.5)--(8.0+1.5,-5.5+1.5);
\node at (8.0+1.5,-4.0+1.5) [fill=red, inner sep=1pt, shape=circle, draw] {};
\node at (8.0+1.5,-4.5+1.5) [fill=red, inner sep=1pt, shape=circle, draw] {};
\node at (8.0+1.5,-5.0+1.5) [fill=red, inner sep=1pt, shape=circle, draw] {};
\node at (8.0+1.5,-5.5+1.5) [fill=red, inner sep=1pt, shape=circle, draw] {};

\draw [->] (9.0,-3.0)-- (12.0,0.0);
\node at (9.0,-3.0)  [shape = circle, draw]{};

\draw[->][red] (8.0,-2.5)--(12,1.5);
\draw[->][red] (8.0,-3.0)--(12,1.0);
\draw[red] (8.0,-2.5)--(8.0,-3.0);
\node at (8.0,-2.5) [fill=red, inner sep=1pt, shape=circle, draw] {};
\node at (8.0,-3.0) [fill=red, inner sep=1pt, shape=circle, draw] {};

\draw[->][blue] (8.0,-3.5)--(12,0.5);
\node at (8.0,-3.5) [fill=red, inner sep=1pt, shape=circle, draw] {};

\draw[->][green] (4.0,-2.10)--(12,5.90);
\draw[->][green] (4.0,-2.60)--(12,5.40);
\draw[->][green] (4.0,-3.10)--(12,4.90);
\draw[->][green] (4.0+3.5,-3.5+3.5)--(12,4.5);
\draw[->][green] (4.0+3.5,-4.10+3.5)--(12,3.90);
\draw[->][green] (4.0+3.5,-4.5+3.5)--(12,3.5);
\draw[->][green] (4.0+3.5,-5.0+3.5)--(12,3.0);
\draw[->][green] (4.0+3.5,-5.5+3.5)--(12,2.5);
\draw[->][green] (4.0+3.5,-6.10+3.5)--(12,1.90);
\draw[->][green] (4.0+3.5,-6.60+3.5)--(12,1.40);
\draw[->][green] (4.0+3.5,-7.10+3.5)--(12,0.90);
\node at (4.0+3.5,-3.5+3.5) [fill=red, inner sep=1pt, shape=circle, draw] {};
\node at (4.0+3.5,-4.10+3.5) [fill=red, inner sep=1pt, shape=circle, draw] {};
\node at (4.0+3.5,-4.5+3.5) [fill=red, inner sep=1pt, shape=circle, draw] {};
\node at (4.0+3.5,-5.0+3.5) [fill=red, inner sep=1pt, shape=circle, draw] {};
\node at (4.0+3.5,-5.5+3.5) [fill=red, inner sep=1pt, shape=circle, draw] {};
\node at (4.0+3.5,-6.10+3.5) [fill=red, inner sep=1pt, shape=circle, draw] {};
\node at (4.0+3.5,-6.60+3.5) [fill=red, inner sep=1pt, shape=circle, draw] {};
\node at (4.0+3.5,-7.10+3.5) [fill=red, inner sep=1pt, shape=circle, draw] {};
\draw[->][red] (4.0+3.5,-3.5+3.5)--(4.0+3.5,-7.10+3.5);

\draw [->] (8.0,-2.0)-- (12.0,2.0);
\node at (8.0,-2.0)  [shape = circle, draw]{};
\draw [->] (7.0,-1.0)-- (12,4.0);
\node at (7.0,-1.0)  [shape = circle, draw]{};
\draw [->] (6.0,0.0)-- (12.0,6.0);
\node at (6.0,0.0)  [shape = circle, draw]{};
\draw [->] (5.0,1.0)-- (12,8.0);
\node at (5.0,1.0)  [shape = circle, draw]{};
\draw [->] (4.0,2.0)-- (12.0,10.0);
\node at (4.0,2.0)  [shape = circle, draw]{};

\draw[->][red] (6.5,0.0)--(12,5.50);
\draw[->][red] (6.5,-0.5)--(12,5.0);
\draw[red] (6.5,0.0)--(6.5,-0.5);
\node at (6.5,0.0) [fill=red, inner sep=1pt, shape=circle, draw] {};
\node at (6.5,-0.5) [fill=red, inner sep=1pt, shape=circle, draw] {};

\draw[->][blue] (4.0,-0.10)--(12,7.90);
\draw[->][blue] (4.0,-0.5)--(12,7.5);
\draw[->][blue] (4.0,-1.0)--(12,7.0);
\draw[->][blue] (4.0,-1.5)--(12,6.5);

\draw[->][red] (4.0,1.5)--(12,9.5);
\draw[->][red] (4.0,1.0)--(12,9.0);
\draw[->][blue] (4.0,0.5)--(12,8.5);

\node at (0.0,-0.25) [fill=yellow, inner sep=1.5pt, shape=diamond, draw] {};
\end{tikzpicture}$$
\caption{The negative block, SS(GNB) and $\mathbb{Z}^{\rm{SS(GNB)}}$ on $BP\mathbb R\langle 3 \rangle^{Q}_{\rost}$: The two curves in orange colours are the exotic multiplications by $a$.} \label{fig:GL2}
\end{figure}
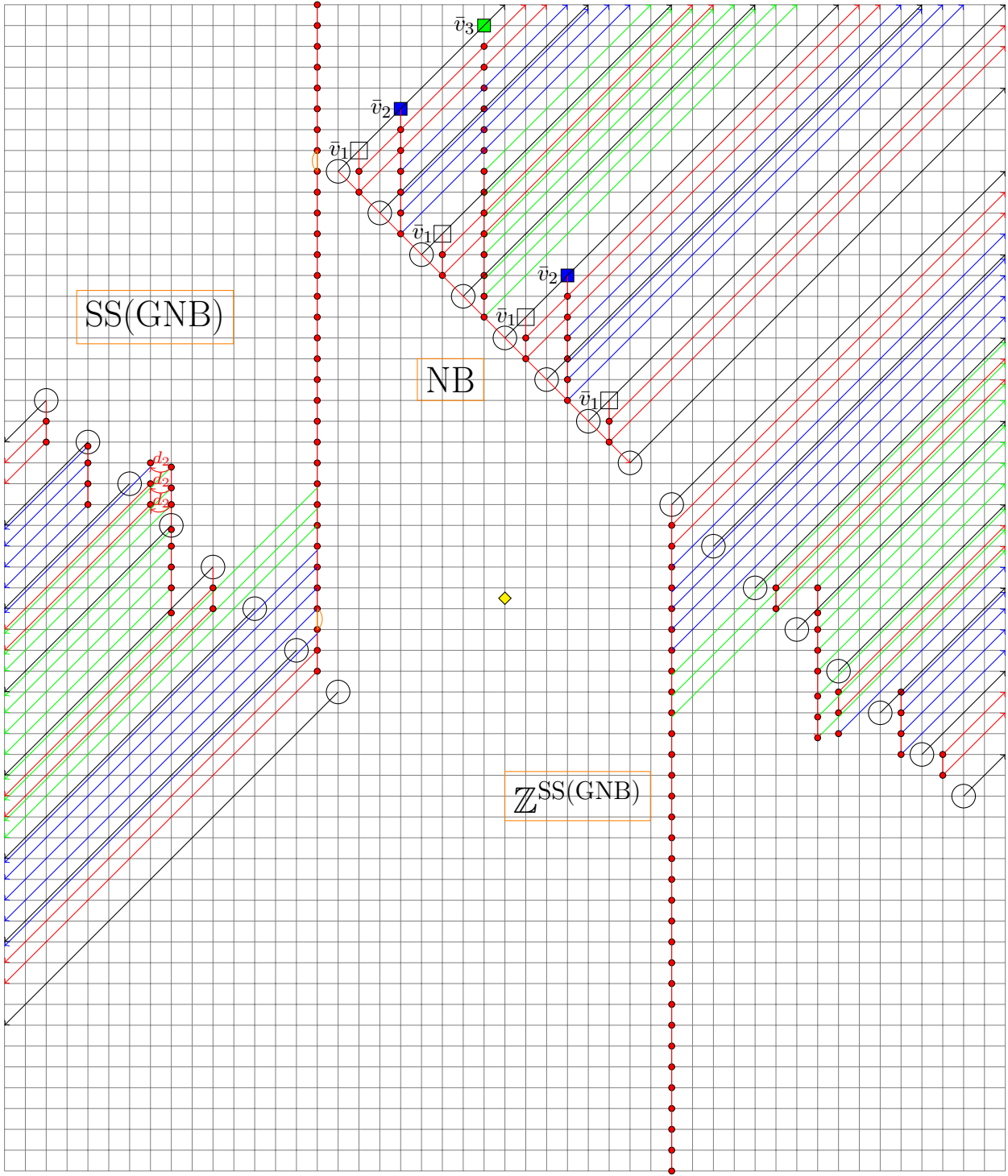

The Figure \ref{fig:GL2} displays the negative block, the spectral sequence of the local cohomology and its Anderson dual. Here, we have mostly
omitted dots, circles and squares except at the ends of diagonals or where an additional generator is required. The vertical lines denote
multiplication by $a$ and the two curved lines in orange colours are the exotic multiplications by $a$ that are not visible on the level of local
cohomology. The yellow diamond does not denote a class, but marks the point we have to reflect (non-torsion classes) at to see Anderson duality. The torsion classes are shifted after reflection by $-1$; that is,
one step horizontally to the left.

The symbol $GNB$ refers to the part of $\Gamma_JBP\R \langle 3 
\rangle$ coming from $NB$, and $SS(GNB)$ is the local cohomology in 
the $E^2$-term that calculates its homotopy. In other words, to go from 
$SS(GNB)$ to $GNB$ there may be differentials and additive 
extensions. 

Once again, we see that the three $d_2$ differentials are forced by
connectivity requirements.

%\newpage

%\newpage
\bigskip

\section{Extension problems}
We give a detailed analysis of the duality as it applies to $BB$.
In other words, we consider each the diagonals $BB_{\delta}$ as
displayed by diagonal in Table \ref{tab:table1}, and calculate its local cohomology
$H^*_J(BB_{\delta})$.  In effect this is $GBB$ (except that  potential
extensions are split, and no differentials are taken into account),
and  these in turn are  displayed by diagonal in Tables
\ref{tab:table2} and \ref{tab:table3}. The reader may wish to annotate this further as a reminder of
the  diagonal $BB_{\delta}$ the local cohomology group came from.

In Table \ref{tab:table2} we have omitted  copies of $\Ftwo$, but we indicate
modules where a differential originates by writing `$d_i$'. In all cases
divisibility of the local cohomology modules shows that once the
differential is known to be non-zero on the top class it is
necessarily a monomorphism on the entire summand. 
When there are several entries in a row in Tables \ref{tab:table2} and \ref{tab:table3}, there is a
possibility of an additive extension. This occurs with
$$\delta \in \{1, 9, 17, 25 \}$$ (the one with $\delta=25$ involves
an undisplayed $\Ftwo$).

The second opportunity for extensions is as part of  Anderson
duality and this also occurs, in our case to make up $NB_8$. The
behaviour on this diagonal is rather
complicated! It involves the diagonals $GBB_{17}$ and $GBB_{18}$.
In fact $GBB_{17}$ already involves a non-split extension
$$0\lra H^3_J(2u^5P)\lra GBB_{17}\lra H^2(\vb_2u^4\Pb_1)/top\lra 0 $$
(where the $/top$ indicates that the differential has killed the top
copy of $\Ftwo$). In other words, we have
$$0\lra P^* \lra GBB_{17}\lra \Sigma^{3\rho} (J\Pb_1)^{\vee}\lra 0, $$
where $$J\Pb_1=(\vb_2, \vb_3)\Pb_1$$ is the augmentation ideal, which
has its bottom class in degree $3\rho$.  The suspensions should be ignored except to say that the top
classes of $P^*$ and $\Sigma^{3\rho} (J \Pb_1)^{\vee}$ are in the
same degree. Actual suspensions are best identified by looking at
Figure \ref{fig:GL1}.
Taking Anderson duals, we have a cofibre sequence of graded abelian groups
$$J\Pb_1\lla P \lla \Z^{GBB_{17}}. $$

The map $$P \lra  J\Pb_1=(\vb_2, \vb_3)$$ in homotopy maps onto the subideal $(\vb_2)$, so that in
fact we have a  splitting
$$\Z^{GBB_{17}}\cong (2, \vb_1)\oplus \Sigma^{4\rho}\Pb_2$$
as abelian groups.  This gives
$$GBB_{17}\cong (2, \vb_1)^*\oplus \Sigma^{-4\rho} \Pb_2^{\vee}. $$
The torsion free part contributes to $$NB_8=(2, \vb_1)$$ and the torsion part to
$$NB_7=(\vb_3)\Pb_2.$$

The  analysis for $NB$ is similar, but a little less complicated.

\newpage

\section{Duality of diagonals}
As described in \cite[Subsection 12.F and Remark 13.3]{BPRn}, one way of thinking about the duality is
that we start with a diagonal $BB_{\delta}$, take local cohomology,
usually in a single degree,  to reach $H^i_J(BB_{\delta})$.  As an
abelian group, this is usually either all torsionfree or annihilated by 2, so
that when Anderson-dualized it contributes to a single diagonal
$NB_{\delta'}$.  We start with $BB_0$, which
gives $NB_{28}$ and then $BB_1, BB_2, \ldots , BB_{28}$ contribute in
succession to $NB_{27}, NB_{26}, \ldots , NB_0$.  Symbolically  we have a close relation of the form
$$H^{3-\epsilon}_J(BB_{\delta})^*\sim NB_{28-\delta -\epsilon'}. $$
However because of differentials, extensions, variation of depth and the question of
torsion, the shift terms $\epsilon, \epsilon'$ are usually non-trivial and we
do not always end up  pairing $\delta, \delta'$ with
$$\delta +\delta'=28.$$ The following table shows what happens, and in
fact 
$$\delta +\delta'\in \{25, 26, 27, 28\}.$$

For example, we see that $\delta =11$ is paired with $\delta'=14$: we
see $$BB_{11}=\Pb_2=\Ftwo [\vb_3]$$ (Table \ref{tab:table3} ($\delta=11$) or Figure \ref{fig:GL1}), with local cohomology
$$H^*_J(BB_{11}) = H^1_J(BB_{11})=\Sigma^{-7\rho}\Pb_2^{\vee}$$ (Table
\ref{tab:table2} ($\delta=10$) or  Figure \ref{fig:GL1}). The Anderson dual of this is $\Pb_2=NB_{14}$ (Table
\ref{tab:table4} ($\delta=14$) or Figure \ref{fig:GL2}). The reader may wish to note why the three
further copies of $\Pb_2$ (in $BB_{12}, BB_{13}$ and $BB_{14}$) do not
lead to copies of $\Pb_2$ in $NB_{13}, NB_{12}$ and $NB_{11}$: in all three
cases there are extensions, in the first case as abelian groups, and
in the others as $\Pb_0$-modules.

\renewcommand{\tabcolsep}{20.00pt}
\renewcommand{\arraystretch}{1.20}
$$\begin{array}{|c|c||c|c|}
\hline
\cellcolor{yellow!55}\delta &\cellcolor{yellow!55}\delta' s.t. H^*_{\Jb_2}(BB_{\delta})^*\sim
NB_{\delta'}&\cellcolor{yellow!55}\delta &\cellcolor{yellow!55}\delta' s.t. H^*_{\Jb_2}(NB_{\delta})^*\sim BB_{\delta'}\\
\hline
0&28&0&28, 25\\
1&26&1&26, 24\\
2&25&2&25, 23\\
3&24&3&24,22 \\
4&24,22&4&24,22,21\\
5&21&5&21, 20\\
6&20&6&20, 19\\
7&d_2&7&d_2\\
8&d_2,20&8&d_2, 20,19\\
9&d_2,18&9&d_2, 18\\
10&16, 17&10&17,16\\
11&14&11&14\\
12&13&12&16,13\\
13&12&13&12\\
14&11&14&11\\
15&d_2&15&\cdot\\
16&d_2, 12,10&16&12, 10\\
17&d_2, 10,9&17&10, 9\\
18&d_2, 9,8&18&9, 8\\
19&d_2, 8&19&7\\
20&d_2, 6,8&20&6, 8\\
21&d_2, 5&21&5\\
22&d_2, 4&22&4\\
23&d_3&23&\cdot\\
24&d_3, 4,3&24&4,3\\
25&2&25&2\\
26&1&26&1\\
27&\cdot&27&\cdot\\
28&0&28&0\\
\hline
\end{array}$$

\bigskip

\section{The last differential} \label{sec:diff}
We have seen that in the local cohomology spectral sequence for
$\BPRn^Q_{\rost}$ there are potential differentials $d_2, \ldots ,
d_n$. Furthermore,  when $n=2$ \cite[Section 13]{BPRn} and $n=3$ each
of these differentials is non-zero somewhere. The purpose of this
section is to show that $d_n\neq 0$ for all
$n\geq 2$.

Focusing on $BB$, let us consider the local cohomology groups
$H^*_J(BB)$. The final entry (i.e.,  the one whose degree has the lowest multiple of
$\sigma$) will be for $2u^{2^n-1}$, which generates a
copy of $P$, and gives rise to $H^n_J(P)\cong P^*$ with top class on
the  $\sigma$ axis at $(n-4(2^n-1))\sigma$. Starting in the next row
up, there is a column of copies of $\Ftwo$ for $-1+y\sigma$ where $-(2^{n+1}-1)\geq
y\geq  n-4(2^n-1)+1$. These entries play an essential role, since the
elements $a^y$ for the same values of $y$ must be cancelled by a differential (this is
necessary by Gorenstein duality since the elements $a^y$ are Anderson dual to entries at $x+y\sigma$ with $x=-(2^{n+1}+1)<0$ and
$x+y<0$).

Since the powers of $a$ are all in 0th local cohomology, to determine
which differential does the cancellation,  we need only 
check which local cohomology group generated the entry on the line $x=-1$. In 
principle we could do this for all of them, but we focus on the
bottom-most class, and give a  little more detail. We
are discussing the entries at
$x-c\sigma$ where $x=0, -1$ and $c=4(2^n-1)-(n+1)=2D_n+n-1$. We already know there is 
$\Ftwo$ from $H^0_J$  (generated by  $a^c$) at the point with $x=0$,  and that there is an
$\Ftwo$ at the point with $x=-1$ coming from $H^n_J(\Pb)$ where $\Pb$ is the principal ideal
generated by $u^{2^n-2}\vb_1a^2$. It remains to
show that this is the entire group at $x=-1$ since then $d_i(a^c)=0$
for $i<n$, and we must have $d_n(a^c)\neq 0$. To see this, we
estimate that the top class 
 of the local cohomology groups coming from
the ideals along the $u$-power line $x+y=0$ (before the local
cohomology shift), must be on or above the line 
$x+y=|\vb_1|+\cdots +|\vb_n|=D_n$. The horizontal distance from the point
at issue (viz $-1-c\sigma$) to the line is precisely $n$ (the
largest possible local cohomological shift), so there is only one
diagonal (viz $\delta=2(2^n-1)-2$) that can contribute to it.

The argument that $a^{c-1}$ also supports a non-zero $d_n$ is precisely
similar.  Lower powers of $a$ support $d_i$ for $i<n$, but we have not
determined the exact values of $i$ for $n\geq 4$.

%\bibliographystyle{alpha}
%\bibliography{Chromatic}

\bigskip

\begin{thebibliography}{100}
\bibitem{Ati66} M.F. Atiyah, {\em $K$-theory and reality},
  Quart. J. Math. Oxford Ser.(2), 17 (1966), 367-386.
\bibitem{AtiyahSegal} M.F. Atiyah and G.B.Segal,  {\em Equivariant
    $K$-theory and completion},
J.Diff.Geom {\bf 3} (1969), 1-18
%\bibitem{Ban13} R. Banerjee, {\em On the $ER(2)$-cohomology of some odd-dimensional projective spa%ces}, Topology Appl.  160(12) (2013), 1395-1405.
\bibitem{BC76} Edgar H. Brown, Jr. and M. Comenetz, {\em Pontrjagin duality for generalized homology and cohomology theories}, Amer. J. Math. 98(1) (1976), 1-27.
\bibitem{BG10} R.R. Bruner and J.P.C. Greenlees, {\em Connective real $K$-theory  of finite groups}, Mathematical Surveys and Monographs, Vol. 169,  American Mathematical Society, Providence, RI, 2010, 318+v pp.
\bibitem{Dug05} D. Dugger, {\em An Atiyah-Hirzebruch spectral sequence for $KR$-theory}, $K$-Theory 35(3-4) (2005), 213-256.
\bibitem{DGI06} W.G. Dwyer, J.P.C. Greenlees and S. Iyengar, {\em Duality in algebra and topology}, Adv. Math.  200(2) (2006), 357-402.
\bibitem{EKMM97} A.D. Elmendorf, I. Kriz, M.A. Mandell and J.P. May, {\em Rings, modules, and algebras in stable homotopy theory}, Mathematical Surveys and Monographs, Vol. 47, American Mathematical Society, Providence, RI, 1997, xii+249 pp.
\bibitem{KEG}
  J.P.C.Greenlees,  
  {\em K-homology of universal spaces and local cohomology of the representation
  ring}
  Topology 
  {\bf 32} (1993) 295-308.
\bibitem{groupca}
  J.P.C.Greenlees,  
  {\em Commutative algebra in group cohomology.}
  J.Pure and Applied Algebra 
  {\bf 98} (1995) 151-162
%\bibitem{G87} J.P.C. Greenlees, {\em Representing Tate cohomology of $G$-spaces}, Proc. Edinburgh Math. Soc. (2) 30(3) (1987), 435-443.
%\bibitem{G88} J.P.C. Greenlees, {\em Stable maps into free $G$-spaces}, Trans. Amer. Math. Soc. 310(1) (1988), 199-215.
\bibitem{hica} J.P.C. Greenlees, {\em Homotopy invariant commutative algebra over fields}, arXiv:1601.02473.
\bibitem{zkr} J.P.C. Greenlees, {\em Four approaches to cohomology theories with reality}, arXiv:1705.09365.
\bibitem{ringlct}
   J.P.C.Greenlees and G.Lyubeznik,  
   {\em Rings with a local cohomology theorem and applications to cohomology 
   rings of groups.}
   J. Pure and Applied Algebra {\bf 149} (2000) 267-285
%\bibitem{GM95a} J.P.C. Greenlees and J.P. May, {\em Completions in algebra and topology}, in Handbook of algebraic topology, North-Holland, Amsterdam, 1995, 255-276.
%\bibitem{GM95b} J.P.C. Greenlees and J.P. May, {\em Generalized Tate cohomology}, Mem. Amer. Math. Soc. 113(543), 1995, viii+178 pp.
\bibitem{BPRn} J.P.C. Greenlees and L. Meier, {\em Gorenstein duality
    for real spectra},  Algebr. Geom. Topol. 17-6 (2017), 3547--3619. DOI 
10.2140/agt.2017.17.3547,
(arXiv:1607.02332).
\bibitem{HHR} M.A. Hill, M.J. Hopkins and D.C. Ravenel, {\em On the nonexistence of elements of Kervaire invariant one}, Ann. of Math. (2) 184(1) (2016), 1-262.
\bibitem{HM} M.A. Hill and L. Meier, {\em The $C_2$-spectrum
    $TMF_1(3)$ and its invertible modules}, 
  Algebr. Geom. Topol.  17 (2017), no. 4, 1953–2011, arXiv:1507.08115
\bibitem{Hu02} P. Hu, {\em On Real-oriented Johnson-Wilson cohomology}, Algebr. Geom. Topol. 2 (2002), 937-947.
\bibitem{HK01} P. Hu and I. Kriz, {\em Real-oriented homotopy theory
    and an analogue of the Adams-Novikov spectral sequence}, Topology,
  40(2) (2001), 317-399.
\bibitem{Karoubi} M.Karoubi,  {\em Sur la $K$-th\'eorie \'equivariante}
  Lecture Notes in Math., {\bf 136} Springer-Verlag (1970) 187-253 
%\bibitem{KW15} N. Kitchloo and W.S. Wilson, {\em The $ER(n)$-cohomology of $BO(q)$  and real Johnso%n-Wilson orientations for vector bundles}, Bull. Lond. Math. Soc.  47(5) (2015), 835-847.
\bibitem{L68} P.S. Landweber, {\em Conjugations on complex manifolds and equivariant homotopy of MU}, Bull. Amer. Math. Soc. 74 (1968), 271-274.
%\bibitem{LO16} G. Laures and M. Olbermann, {\em $TMF_0 (3)$-characteristic classes for string bundles}, Math. Z. 282(1-2) (2016), 511-533.
%\bibitem{Lor16} V. Lorman, {\em The real Johnson-Wilson cohomology of $\mathbb{C}\mathbb{P}^\infty$}, Topology Appl.  209 (2016), 367-388.
%\bibitem{Ric16} N. Ricka, {\em Equivariant Anderson duality and Mackey functor duality},
  %  Glasg. Math. J.  58(3) (2016), 649-676.

\end{thebibliography}
\end{document}